\documentclass{amsart}

\title{Pre-model Structures on the Category of \texorpdfstring{$C_p$}{Cₚ}-Sets}
\author[D. Mehrle]{David Mehrle}
\address{Colorado State University, Fort Collins, CO, USA}
\email{david.mehrle@colostate.edu}

\author[B. Spitz]{Ben Spitz}
\address{Indiana University, Bloomington, IN, USA}
\email{bespitz@iu.edu}


\usepackage{amsmath,amssymb,amsthm}

\usepackage[colorlinks=true,linkcolor=frenchblue,citecolor=forestgreen]{hyperref}
\usepackage[nameinlink,capitalize]{cleveref}

\usepackage{quiver}
\usepackage[shortlabels]{enumitem}
\usepackage{stmaryrd}
\usepackage{mathtools}
\usepackage{xcolor}
\usepackage{microtype}
\usepackage[margin=1in]{geometry}

\usepackage{tabularray}
\UseTblrLibrary{booktabs}

\usepackage{rotating}
\usepackage{fontawesome5}
\newcommand{\computerassisted}{\faLaptopCode}
\usepackage{standalone}

\definecolor{forestgreen}{rgb}{0.13, 0.55, 0.13}
\definecolor{frenchblue}{rgb}{0.0, 0.45, 0.73}
\definecolor{firered}{rgb}{0.86,0.00,0.00}

\newcommand{\cmark}{\checkmark}
\newcommand{\xmark}{\text{\sffamily X}}
\renewcommand\emptyset\varnothing

\newcommand{\Set}{\mathsf{Set}}
\newcommand{\GSet}{\Set^G}
\newcommand{\CpSet}{\Set^{C_p}}
\newcommand{\cat}{\mathcal}
\newcommand{\catC}{\mathcal{C}}
\newcommand{\catD}{\mathcal{D}}

\newcommand{\llp}[1]{\prescript{\boxslash}{}{#1}}
\newcommand{\rlp}[1]{{#1}^{\boxslash}}

\newcommand{\all}{\mathsf{all}}
\newcommand{\inj}{\mathsf{inj}}
\newcommand{\surj}{\mathsf{surj}}
\newcommand{\bij}{\mathsf{bij}}
\newcommand{\init}{\mathsf{init}}
\newcommand{\pb}{\mathsf{pb}}
\newcommand{\homog}{\mathsf{homog}}

\DeclareMathOperator{\id}{id}

\crefname{equation}{}{}

\numberwithin{equation}{section}
\newtheorem{theorem}[equation]{Theorem}
\AddToHook{env/theorem/begin}{\crefalias{equation}{theorem}}
\newtheorem*{theorem*}{Theorem}
\newtheorem{proposition}[equation]{Proposition}
\AddToHook{env/proposition/begin}{\crefalias{equation}{proposition}}
\newtheorem{corollary}[equation]{Corollary}
\AddToHook{env/corollary/begin}{\crefalias{equation}{corollary}}

\AddToHook{env/conjecture/begin}{\crefalias{equation}{conjecture}}
\newtheorem{lemma}[equation]{Lemma}
\AddToHook{env/lemma/begin}{\crefalias{equation}{lemma}}
\theoremstyle{definition}
\newtheorem{definition}[equation]{Definition}
\AddToHook{env/definition/begin}{\crefalias{equation}{definition}}
\newtheorem{remark}[equation]{Remark}
\AddToHook{env/remark/begin}{\crefalias{equation}{remark}}
\newtheorem{example}[equation]{Example}
\AddToHook{env/example/begin}{\crefalias{equation}{example}}
\newtheorem{strategy}[equation]{Strategy}
\crefname{strategy}{Strategy}{Strategies}
\AddToHook{env/strategy/begin}{\crefalias{equation}{strategy}}

\newcommand*{\IntersectionTable}[1]{%
    \begin{tblr}{
        colspec={l|*{6}{l}},
        colsep=2pt,
        cells={mode=math},
        #1
    }
		\toprule
		\downarrow \cap \rightarrow                          & \inj^e                                               & \inj^{C_p}                                           & \surj^e                                & \surj^{C_p}                            & \pb                                                  & \homog                                               \\ \midrule
		\init^{C_p}                                          & \init^{C_p} \cap \inj^e                              & \init^{C_p}                                          & \init^{C_p} \cap \surj^e               & \{\id_\emptyset\}^{C_p}              & \init^{C_p} \cap \pb                                 & \init^{C_p}                                          \\
		\init^{C_p} \cap \pb                                 & \init^{C_p} \cap \inj^e                              & \init^{C_p} \cap \pb                                 & \{\id_\emptyset\}^{C_p} \cap \surj^e & \{\id_\emptyset\}^{C_p}              & \init^{C_p} \cap \pb                                 & \init^{C_p} \cap \pb                                 \\
		\init^e \cup \{\id_\emptyset\}^{C_p}               & (\init^e \cup \{\id_\emptyset\}^{C_p}) \cap \inj^e & \init^e \cup \{\id_\emptyset\}^{C_p}               & \{\id_\emptyset\}^{C_p} \cap \surj^e & \{\id_\emptyset\}^{C_p}              & \init^e \cup \{\id_\emptyset\}^{C_p}               & \init^e \cup \{\id_\emptyset\}^{C_p}               \\
		\init^e                                              & \init^e                                              & \init^e                                              & \{\id_\emptyset\}                    & \{\id_\emptyset\}^{C_p} \cap \init^e & \init^e                                              & \init^e                                              \\
		\{\id_\emptyset\}^{C_p}                            & \{\id_\emptyset\}^{C_p} \cap \inj^e                & \{\id_\emptyset\}^{C_p}                            & \{\id_\emptyset\}^{C_p} \cap \surj^e & \{\id_\emptyset\}^{C_p}              & \{\id_\emptyset\}^{C_p}                            & \{\id_\emptyset\}^{C_p}                            \\
		\{\id_\emptyset\}^{C_p} \cap \init^e               & \{\id_\emptyset\}^{C_p} \cap \init^e               & \{\id_\emptyset\}^{C_p} \cap \init^e               & \{\id_\emptyset\}                    & \{\id_\emptyset\}^{C_p} \cap \init^e & \{\id_\emptyset\}^{C_p} \cap \init^e               & \{\id_\emptyset\}^{C_p} \cap \init^e               \\
		\{\id_\emptyset\}                                  & \{\id_\emptyset\}                                  & \{\id_\emptyset\}                                  & \{\id_\emptyset\}                    & \{\id_\emptyset\}                    & \{\id_\emptyset\}                                  & \{\id_\emptyset\}                                  \\
		\init^{C_p} \cap \inj^e                              & \init^{C_p} \cap \inj^e                              & \init^{C_p} \cap \inj^e                              & \{\id_\emptyset\}^{C_p} \cap \bij^e  & \{\id_\emptyset\}^{C_p} \cap \inj^e  & \init^{C_p} \cap \inj^e                              & \init^{C_p} \cap \inj^e                              \\
		(\init^e \cup \{\id_\emptyset\}^{C_p}) \cap \inj^e & (\init^e \cup \{\id_\emptyset\}^{C_p}) \cap \inj^e & (\init^e \cup \{\id_\emptyset\}^{C_p}) \cap \inj^e & \{\id_\emptyset\}^{C_p} \cap \bij^e  & \{\id_\emptyset\}^{C_p} \cap \inj^e  & (\init^e \cup \{\id_\emptyset\}^{C_p}) \cap \inj^e & (\init^e \cup \{\id_\emptyset\}^{C_p}) \cap \inj^e \\
		\init^{C_p} \cap \surj^e                             & \{\id_\emptyset\}^{C_p} \cap \bij^e                & \init^{C_p} \cap \surj^e                             & \init^{C_p} \cap \surj^e               & \{\id_\emptyset\}^{C_p} \cap \surj^e & \{\id_\emptyset\}^{C_p} \cap \surj^e               & \init^{C_p} \cap \surj^e                             \\
		\{\id_\emptyset\}^{C_p} \cap \inj^e                & \{\id_\emptyset\}^{C_p} \cap \inj^e                & \{\id_\emptyset\}^{C_p} \cap \inj^e                & \{\id_\emptyset\}^{C_p} \cap \bij^e  & \{\id_\emptyset\}^{C_p} \cap \inj^e  & \{\id_\emptyset\}^{C_p} \cap \inj^e                & \{\id_\emptyset\}^{C_p} \cap \inj^e                \\
		\{\id_\emptyset\}^{C_p} \cap \surj^e               & {\{\id_\emptyset\}^{C_p} \cap \bij^e}              & \{\id_\emptyset\}^{C_p} \cap \surj^e               & \{\id_\emptyset\}^{C_p} \cap \surj^e & \{\id_\emptyset\}^{C_p} \cap \surj^e & \{\id_\emptyset\}^{C_p} \cap \surj^e               & \{\id_\emptyset\}^{C_p} \cap \surj^e               \\
		\{\id_\emptyset\}^{C_p} \cap \bij^e                & {\{\id_\emptyset\}^{C_p} \cap \bij^e}              & \{\id_\emptyset\}^{C_p} \cap \bij^e                & \{\id_\emptyset\}^{C_p} \cap \bij^e  & \{\id_\emptyset\}^{C_p} \cap \bij^e  & \{\id_\emptyset\}^{C_p} \cap \bij^e                & \{\id_\emptyset\}^{C_p} \cap \bij^e                \\
		\bottomrule
    \end{tblr}%
}

\begin{document}

\begin{abstract}
	It is a folk theorem of Goodwillie that there are exactly 9 model structures on the category of sets, coming from 18 pre-model structures. 
	We investigate model structures and pre-model structures on the category of $G$-sets for a finite group $G$. 
	We outline a strategy for finding all model structures in this category, and partially execute this strategy  for the case of $C_p$-sets, where $C_p$ denotes a cyclic group of prime order $p$. 
	In particular, there are exactly 1777 pre-model structures on the category of $C_p$-sets.
\end{abstract}

\maketitle

\setcounter{tocdepth}{1}
\tableofcontents

\vspace*{-2em}

\section{Introduction}

Model structures were introduced by Quillen \cite{Quillen1967} to axiomatize homotopy theory. Although Quillen was interested in a uniform treatment of familiar homotopy theories, the field of homotopical combinatorics has found reason to study model structures as combinatorial objects. Previous work studies model structures on finite total orders \cite{BOOR2023}, finite lattices \cite{MORSVZ2025}, posets \cite{DZ2021},  graphs \cite{Droz2012}, directed graphs \cite{CDKOSW2022} and the category of sets \cite{ACB2025}. Surprisingly, a total classification of model structures in the category of sets is possible: it is a result of Goodwillie that there are exactly nine \cite{ACB2025}.

We are interested in model structures on the category of $G$-sets and equivariant functions, where $G$ is a finite group. Unfortunately, the task of classifying model structures is computationally infeasible for arbitrary groups. Nevertheless, we make partial progress towards the count of model structures for a group $G = C_p$ of prime order by counting the pre-model structures on the category of $C_p$-sets. The~answer~doesn't~depend~on~$p$.

\begin{theorem*}[\cref{cor:final count}]
	Let $p$ be a prime number and let $C_p$ be a cyclic group of order $p$. There are exactly $1777$ pre-model structures on the category of $C_p$-sets, independent of $p$. 
\end{theorem*}

To arrive at this number, we compute the lattice of weak factorization systems on the category of $C_p$-sets in \Cref{fig:biggest lattice} and count the intervals in this lattice by computer. Each interval corresponds to a pre-model structure.
Our code is available on GitHub \cite{code}. Results whose proofs rely on computer verification are marked with the symbol \computerassisted.

By contrast, there are only $18$ pre-model structures on the category of sets and $9$ pre-model structures on the category of $k$-vector spaces (independent of the field $k$). These can be found by counting the intervals in the posets of weak factorization systems for these categories \cite[Diagrams 3.5 and 6.5]{ACB2025}.

We leave open the question of which of these pre-model structures are model structures, and leave open the description of the fibrations and cofibrations in these model structures.

\subsection*{Acknowledgments}

The authors would like to thank Valentina Zapata Castro for asking this question, David Chan for his enthusiasm and interest, and Omar Antol\'in-Camarena and Tobias Barthel for encouraging us to finish this project. We also thank an anonymous referee for their many insightful and detailed comments which greatly improved not only this paper, but also the accompanying code.

\subsection*{AI Use}

The authors of this article used AI for neither the mathematics nor the writing presented here. However, we did use AI to help us write the accompanying computer code \cite{code}. We used AI to help with code generation and debugging, but the algorithms presented there are our own.

\section{A Strategy for Classifying Model Structures}
\label{section: the strategy}

To prove this theorem, we follow a strategy similar to that outlined in \cite{ACB2025}, which relies on a definition of model structures in terms of weak factorization systems. We first recall the definition of lifting properties, which are inseparable from model structures.

\begin{definition}
	Let $\catC$ be any category. We say that a pair of morphisms $(f,g)$ in $\catC$ has the \emph{lifting property} if, for all solid commutative squares
	
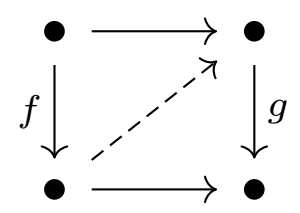
\[\begin{tikzcd}[ampersand replacement=\&]
			\bullet \& \bullet \\
			\bullet \& \bullet
			\arrow[from=1-1, to=1-2]
			\arrow["f"', from=1-1, to=2-1]
			\arrow["g", from=1-2, to=2-2]
			\arrow[dashed, from=2-1, to=1-2]
			\arrow[from=2-1, to=2-2]
		\end{tikzcd}\]
	there exists a dashed lift making the two triangles commute. If this is the case, we write $f \boxslash g$ and say that $f$ has the \emph{left lifting property} with respect to $g$, or say that $g$ has the \emph{right lifting property} with respect to $f$.
\end{definition}

We are interested in characterizing those morphisms which lift with respect to another class of morphisms.

\begin{definition}
	\label{def of lifting classes}
	Let $\catC$ be a category, and let $\mathcal{M}$ be a class of morphisms in $\catC$.
	\begin{enumerate}[(a)]
		\item Let $\rlp{\mathcal M}$ be the \emph{right-lifting class of $\mathcal{M}$}: the class of morphisms $g$ such that $g$ has the right lifting property with respect to all morphisms in $\mathcal{M}$.
		      \[
			      \rlp{\mathcal{M}} := \{g \in \catC \mid f \boxslash g \text{ for all } f \in \mathcal M\}
		      \]
		\item Let $\llp{\mathcal M}$ be the \emph{left-lifting class of $\mathcal{M}$}: the class of morphisms $f$ such that $f$ has the left lifting property with respect to all morphisms in $\mathcal M$:
		      \[
			      \llp{\mathcal{M}} := \{f \in \catC \mid f \boxslash g \text{ for all } g \in \mathcal M\}.
		      \]
	\end{enumerate}
\end{definition}

\begin{example}
	\label{lifts of folds}
	Let $I$ be a (possibly infinite) set, and let $n = \lvert I \rvert$. Consider the (unique) morphism $\nabla \colon I \to \ast$ in $\Set$, where $\ast$ denotes a singleton set. 
	\begin{itemize}
		\item If $n = 0$, then $\rlp{\{\nabla\}} = \rlp{\big\{\emptyset \to \ast\big\}}$ is the class $\surj$ of all surjections in $\Set$, because diagrams
		      
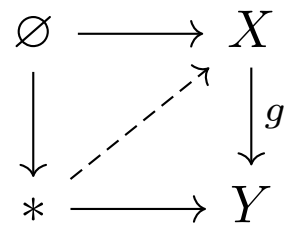
\[
			      \begin{tikzcd}
				      \emptyset \ar[r] \ar[d] & {X} \ar[d, "g"] \\
				      \ast \ar[r] \ar[ur, dashed] & {Y}
			      \end{tikzcd}
		      \]
		      correspond bijectively with elements of $Y$, and lifts correspond bijectively with preimages of those elements under $g$.

		\item If $n = 1$, then $\rlp{\{\nabla\}}$ is the class $\all$ of all morphisms in $\Set$, because $\nabla$ is an isomorphism in this case.

		\item If $n > 1$, then $\rlp{\{\nabla\}}$ is the class $\inj$ of all injections in $\Set$, because diagrams
		      
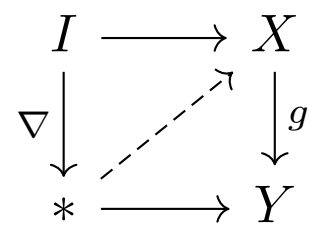
\[
			      \begin{tikzcd}
				      {I} \ar[d, "\nabla"'] \ar[r] & X \ar[d, "g"] \\
				      \ast \ar[r] \ar[ur, dashed] & Y
			      \end{tikzcd}
		      \]
		      correspond bijectively with tuples $(x_i)_{i \in I}$ of elements of $X$ with common image under $g$, for which a lift exists if and only if all elements in the tuple are equal.
	\end{itemize}
\end{example}

A pair of classes of morphisms that lift against each other is a weak factorization system if any morphism can be factored into a composite of morphisms, one in each class. More precisely:

\begin{definition}
	A \emph{weak factorization system} on $\catC$ is a pair $(\cat L, \cat R)$ of classes of morphisms in $\catC$ such that
	\begin{enumerate}[(a)]
		\item every morphism in $\catC$ can be factored as $gf$ with $f \in \cat L$ and $g \in \cat R$, and
		\item $\cat L$ is the class of morphisms with the left-lifting property against $\cat R$, and $\cat R$ is the class of morphisms with the right-lifting property against $\cat L$:
		      \[
			      \cat L = \llp{\cat R}
			      \quad \text{ and } \quad
			      \rlp{\cat L} = \cat R.
		      \]
	\end{enumerate}
\end{definition}

Although model structures are typically presented in terms of fibrations, cofibrations, and weak equivalences, the following definition in terms of weak factorization systems is better suited to our purposes. See \cite[Definition 7.7]{jt2007}, \cite[Definition 4.1]{ACB2025}, or \cite[Section 2]{BBO2024} for a discussion of how this definition of a model category coincides with the standard definition.

\begin{definition}\
	\begin{enumerate}[(a)]
		\item A \emph{pre-model structure} on a category $\catC$ is a pair of weak factorization systems $(\cat L, \cat R)$ and $(\cat L', \cat R')$ such that $\cat R \subseteq \cat R'$ (or equivalently, $\cat L \supseteq \cat L'$).
		\item A pre-model structure is a \emph{model structure} if additionally, the class $\cat{W} := \cat R \circ \cat L'$ satisfies the two-out-of-three property: given two composable morphisms $g$ and $f$, if any two of $f, g$, and $g \circ f$ are in $\mathcal{W}$, then so is the third.
	\end{enumerate}
\end{definition}

In light of this definition, the classification of model structures requires us to classify the pre-model structures, which in turn requires a classification of the weak factorization systems. A key simplification in the classification of weak factorization systems comes from the following incarnation of the small object argument. The small object argument is originally due to Quillen \cite[Lemma 2.3.4]{Quillen1967}; we recall a reformulation from \cite[Theorem 12.2.2 and Corollary 12.2.4]{Riehl2014}.

\begin{theorem}[Small Object Argument]
	\label{smallobjectargument}
	Let $\catC$ be a locally presentable category, and let $I$ be a \emph{set} of morphisms in $\catC$. Then
	\[
		\left(\llp{(\rlp{I})},\rlp{I}\right)
	\]
	is a weak factorization system on $\catC$.
\end{theorem}

\begin{remark}
	Note that $\GSet$ is locally presentable (it is the category of presheaves on the small category $BG$), so the small object argument as stated above applies to $\GSet$.
\end{remark}

Any weak factorization system $(\cat L, \cat R)$ is superficially of this form (take $I = \cat L$); however, $\cat L$ is typically a proper class rather than a set. However, if it turns out that every right-lifting class $\rlp{\mathcal M}$ of a class of morphisms $\mathcal M$ in $\catC$ is also the right-lifting class $\rlp{I}$ of some \emph{set} $I$ of morphisms, then classifying weak factorization systems on $\catC$ is equivalent to classifying right-lifting classes, via $(\cat L, \cat R) \mapsto \cat R$. We will see that this is the case for $\CpSet$ (\cref{thm:wfs lattice}).
By \cref{lifting class of union} below, the right-lifting class of a class of morphisms $I$ is equivalently described as the intersection of right-lifting classes with respect to each individual morphism in $I$:
\[
	\rlp{I} = \bigcap_{f \in I} \rlp{\{f\}}.
\]
So a sensible place to start is to find all right-lifting classes against individual morphisms. In general, finding all possible right-lifting classes is a formidable task.

\section{Lifting Properties}

We record a few properties of lifting classes below.

\begin{proposition}\label{proposition:all lift against isos}
	If $f$ is an isomorphism, then
	\[
		\rlp{\{f\}} = \llp{\{f\}} = \operatorname{Mor}(\catC).
	\]
\end{proposition}

\begin{proposition}
	If $f$ and $f'$ are isomorphic (in the arrow category of $\catC$), then
	\[
		\rlp{\{f\}} = \rlp{\{f'\}}.
	\]
\end{proposition}

\begin{proposition}
	\label{lifting class of union}
	For any class $\mathbf{M}$ of classes of morphisms in a category $\catC$,
	\[
		\rlp{\left(\bigcup_{\mathcal{M} \in \mathbf{M}} \mathcal{M}\right)}
		= \bigcap_{\mathcal{M} \in \mathbf{M}} \rlp{\mathcal{M}}
		\qquad
		\text{ and }
		\qquad
		\llp{\left(\bigcup_{\mathcal{M} \in \mathbf{M}} \mathcal{M}\right)}
		= \bigcap_{\mathcal{M} \in \mathbf{M}} \llp{\mathcal{M}}.
	\]
\end{proposition}

Note that sometimes a morphism lifts against another simply because there are no squares to lift against. In other words, the lifting property is vacuously true.

\begin{definition}
	Let $\mathcal{M}^\emptyset$ be the \emph{vacuous right-lifting class of $\mathcal{M}$}: the class of morphisms $g$ such that, for all morphisms $f$ in $\mathcal{M}$, there does not exist any commutative square with $f$ on the left and $g$ on the right:
	\[
		\mathcal{M}^\emptyset :=
		\left\{
		g \in \catC
		\ \middle\vert\
		\text{for all } f \in \mathcal{M} \text{, there are no commutative squares }
		\begin{tikzcd}[cramped, sep=1em]
			\bullet \ar[d, "f"'] \ar[r] & \bullet \ar[d, "g"] \\
			\bullet \ar[r] & \bullet
		\end{tikzcd}
		\right\}.
	\]
	For emphasis: one always has $\mathcal{M}^\emptyset \subseteq \mathcal{M}^\boxslash$, and usually this containment is proper.
\end{definition}

\begin{example}
	Consider the class $\mathcal M = \textsf{Id}$ of all identity morphisms in $\Set$. The vacuous right-lifting class $\mathcal M^\emptyset$ is empty in this case; for any $g \colon X \to Y$, there is a commuting diagram
	\[
		\begin{tikzcd}
			\emptyset \ar[r] \ar[d,"\id"'] & X \ar[d, "g"] \\
			\emptyset \ar[r] & Y.
		\end{tikzcd}
	\]
	On the other hand, for an arbitrary class $\mathcal M$ of morphisms in $\Set$, $\id_\emptyset \in \mathcal{M}^\emptyset$ if and only if every morphism in $\mathcal M$ has nonempty codomain. Because the initial object in $\Set$ is strict, the only morphism with empty codomain is $\id_\emptyset$, so this is equivalent to $\id_\emptyset \notin \mathcal M$.
\end{example}

\begin{example}
	Let $G$ be a finite group and let $\GSet$ be the category of $G$-sets and equivariant functions. Let $X$ be a $G$-set with no fixed points. In this category, there are no morphisms $\ast \to X$ because no such function can be equivariant. Consider the class of morphisms $\mathcal M = \{\id_\ast\}$. Then $\mathcal M^\emptyset$ is the class of morphisms $g \colon X \to Y$ such that $X$ has no fixed points; in this case, there are simply no commuting squares
	\[
		\begin{tikzcd}
			\ast \ar[r]\ar[d, "\id"'] & X \ar[d, "g"] \\
			\ast \ar[r] & Y.
		\end{tikzcd}
	\]
\end{example}

The vacuous lifting classes can be used to describe how lifting classes interact with coproducts.

\begin{proposition}
	\label{lifting class of coproduct}
	Suppose $\catC$ admits coproducts. For any set $\{f_i \mid i \in I\}$ of morphisms in $\catC$, we have:
	\begin{enumerate}[(a)]
		\item $\displaystyle \left\{\coprod_{i \in I} f_i\right\}^\emptyset = \bigcup_{i \in I} \{f_i\}^\emptyset$;
		\item $\displaystyle
			      \rlp{\left\{\coprod_{i \in I} f_i\right\}} = \rlp{\{f_i \mid i \in I\}} \cup \bigcup_{i \in I} \{f_i\}^\emptyset.
		      $
	\end{enumerate}
\end{proposition}

\begin{proof}
	First, note that commutative diagrams
	\[\begin{tikzcd}[ampersand replacement=\&]
			\bullet \& \bullet \\
			\bullet \& \bullet
			\arrow[from=1-1, to=1-2]
			\arrow["{\displaystyle \coprod_{i\in I} f_i}"', from=1-1, to=2-1]
			\arrow["g", from=1-2, to=2-2]
			\arrow[dashed, from=2-1, to=1-2]
			\arrow[from=2-1, to=2-2]
		\end{tikzcd}\]
	correspond bijectively with tuples of commutative diagrams
	\begin{equation}
		\label{tuples of squares}
		\left(
		\begin{tikzcd}[ampersand replacement=\&, cramped, sep=5mm]
		\bullet \& \bullet \\
		\bullet \& \bullet
		\arrow[from=1-1, to=1-2]
		\arrow["{f_i}"', from=1-1, to=2-1]
		\arrow["g", from=1-2, to=2-2]
		\arrow[dashed, from=2-1, to=1-2]
		\arrow[from=2-1, to=2-2]
		\end{tikzcd}
		\right)_{\!\!\!i \in I.}
	\end{equation}
	Claim (a) is then clear -- there is at least one such tuple if and only if $g \notin \{f_i\}^\emptyset$ for all $i$.

	To prove (b), we will first show that $\{f_i \mid i \in I\}^\boxslash \subseteq \{\coprod_i f_i\}^\boxslash$ and $\bigcup_{i \in I} \{f_i\}^\emptyset \subseteq \{\coprod_i f_i\}^\boxslash$.

	Assume that $g \in \rlp{\{f_i \mid i \in I\}}$. Then $f_i \boxslash g$ for all $i \in I$, so by the bijection above $\coprod_i f_i \boxslash g$ as well. Therefore,
		\[
			\rlp{\{f_i \mid i \in I\}} \subseteq \rlp{\left\{\coprod_{i \in I} f_i\right\}}.
		\]

Now by (a),
	\begin{equation}
		\label{union of vacuous is coprod vacuous}
		\bigcup_{i \in I} \{f_i\}^\emptyset = \left\{\coprod_{i \in I} f_i\right\}^\emptyset \subseteq \rlp{\left\{\coprod_{i \in I} f_i\right\}}.
	\end{equation}
	This establishes the containment from right-to-left.
	
	\medskip

	In the other direction, assume that $g \in \rlp{\{\coprod_i f_i\}}$. There are two cases: either $g \in \{\coprod_i f_i\}^\emptyset$ or $g \not \in \{\coprod_i f_i\}^\emptyset$.

	In the first case, by (a),
	\begin{equation}
		\label{coprod vacuous is union of individual vacuous}
		g \in \left\{\coprod_{i \in I} f_i\right\}^\emptyset = \bigcup_{i \in I} \{f_i\}^\emptyset.
	\end{equation}

	In the second case, we assume $g \in \rlp{\{\coprod_i f_i\}}$ and $g \not \in \{\coprod_i f_i\}^\emptyset$. By (a), we have $g \notin \{f_i\}^{\emptyset}$ for all $i$. Fix for each $i$ a square $S_i$ witnessing that $g \notin \{f_i\}^\emptyset$.
	
	Now fix some $j \in I$, and suppose we're given a single square
	\[\begin{tikzcd}[ampersand replacement=\&]
			\bullet \& \bullet \\
			\bullet \& \bullet.
			\arrow[from=1-1, to=1-2]
			\arrow["{f_j}"', from=1-1, to=2-1]
			\arrow["g", from=1-2, to=2-2]
			\arrow[from=2-1, to=2-2]
			\arrow["{S'}"{description}, draw=none, from=1-1, to=2-2]
		\end{tikzcd}\]
	By the previous paragraph, there also exist squares
	\[\begin{tikzcd}
			\bullet & \bullet \\
			\bullet & \bullet
			\arrow[from=1-1, to=1-2]
			\arrow[""{name=0, anchor=center, inner sep=0}, "{f_i}"', from=1-1, to=2-1]
			\arrow[""{name=1, anchor=center, inner sep=0}, "g", from=1-2, to=2-2]
			\arrow[from=2-1, to=2-2]
			\arrow["{S_i}"{description}, draw=none, from=0, to=1]
		\end{tikzcd}\]
	for all $i \in I$ with $i \neq j$. Combining these yields a tuple of squares as in \cref{tuples of squares}, and therefore a square
	\[
		\begin{tikzcd}
			\bullet \ar[d, "\displaystyle\coprod_{i \in I} f_i"'] \ar[r] & \bullet \ar[d, "g"] \\
			\bullet \ar[r] & \bullet
		\end{tikzcd}
	\]
	which admits a lift by assumption. Passing through the bijection once again, this furnishes each square
	$\begin{tikzcd}[sep=10, cramped]
			\bullet \ar[r] \ar[d, "f_i"'] & \bullet \ar[d, "g"] \\
			\bullet \ar[r] & \bullet
		\end{tikzcd}
	$
	with a lift. In particular, $S'$ admits a lift. Since $S'$ was arbitrary, $f_j \boxslash g$. Hence, $g \in \rlp{\{f_j\}}$.

	The same argument works for any $j \in I$, so we conclude that $g \in \rlp{\{f_j\}}$ for all $j \in I$.
	\[
		g \in \{f_i \mid i \in I\}^\boxslash
	\]
	by an application of \cref{lifting class of union}.
\end{proof}

Lifting classes can be carried across adjunctions.

\begin{proposition}\label{proposition:adjunction boxslash}
	Let $L \dashv R$ be an adjoint pair of functors, with $L : \catC \to \catD$ and $R : \catD \to \catC$. Let $\mathcal{M}$ be a class of morphisms in $\catC$. Then
	\[
		\rlp{L(\mathcal{M})} = R^{-1}(\rlp{\mathcal{M}})
		\quad \text{ and } \quad
		L(\mathcal{M})^\emptyset = R^{-1}(\mathcal{M}^\emptyset).
	\]
\end{proposition}

\begin{proof}
	Commutative diagrams
	in $\catC$ as on the left correspond bijectively with commutative diagrams in $\cat D$ as on the right.
	\[
		\begin{tikzcd}[ampersand replacement=\&]
			\bullet \& \bullet \\
			\bullet \& \bullet
			\arrow[from=1-1, to=1-2]
			\arrow["f"', from=1-1, to=2-1]
			\arrow["{R(g)}", from=1-2, to=2-2]
			\arrow[dashed, from=2-1, to=1-2]
			\arrow[from=2-1, to=2-2]
		\end{tikzcd}
		\qquad
		\begin{tikzcd}[ampersand replacement=\&]
			\bullet \& \bullet \\
			\bullet \& \bullet
			\arrow[from=1-1, to=1-2]
			\arrow["{L(f)}"', from=1-1, to=2-1]
			\arrow["g", from=1-2, to=2-2]
			\arrow[dashed, from=2-1, to=1-2]
			\arrow[from=2-1, to=2-2]
		\end{tikzcd}\qedhere
	\]
\end{proof}

\section{The Category of \texorpdfstring{$G$}{G}-Sets}

Let $G$ be a finite group. We are interested in classifying pre-model structures on $\GSet$, the category of $G$-sets with $G$-equivariant functions between them. Here we review the basic structure of this category, culminating with a structure theorem for morphisms which we later apply to the classification of pre-model structures. For a more thorough review of $G$-sets, we recommend \cite[Chapter 4]{dummitfoote} or \cite[Chapter 1]{tD1979}.

\begin{definition}\
	\begin{enumerate}[(a)]
		\item A (left) \emph{$G$-set} is a set $X$ equipped with an associative and unital function $G \times X \to X$, written $(g, x) \mapsto g \cdot x$. This function is called a \emph{$G$-action}.
		\item A \emph{morphism of $G$-sets} $f \colon X \to Y$ is an equivariant function: a function that commutes with the actions on the domain and codomain.
	\end{enumerate}
	Let $\GSet$ be the category of (left) $G$-sets and equivariant functions.
\end{definition}

\begin{example}
	Let $H$ be a subgroup of $G$. The set of left cosets
	\[
		G/H := \{gH \mid g \in G\}
	\]
	is a $G$-set, with $G$-action by $a \cdot gH = (ag)H$.
\end{example}

This example is fundamental in the sense that every $G$-set is isomorphic to a disjoint union of sets of left cosets, as we will see below.

\begin{definition}
	Let $X$ be a $G$-set with an element $x\in X$.
	\begin{enumerate}[(a)]
		\item The \emph{orbit} of $x$ is the sub-$G$-set
		      \[
			      G \cdot x := \{g \cdot x \mid g \in G\} \subseteq X.
		      \]
		\item The \emph{stabilizer} of $x$ is the subgroup
		      \[
			      G_x := \{g \in G \mid g \cdot x = x\} \leq G.
		      \]
	\end{enumerate}
\end{definition}

Orbits and stabilizers are related by the following.

\begin{proposition}[Orbit-Stabilizer Theorem]
	Given a $G$-set $X$ and an element $x \in X$, there is an isomorphism of $G$-sets between the orbit $G \cdot x$ of $x$ and the coset space $G/G_x$.
\end{proposition}

\begin{proposition}
	Let $X$ be a $G$-set. The set of orbits of elements of $X$ forms a partition of $X$, i.e. two orbits are either disjoint or equal.
\end{proposition}

Consequently, any $G$-set $X$ can be written as a disjoint union of $G$-sets isomorphic to $G/H$. We refer to the terms of this disjoint union as the \emph{orbits of $X$}. It is also common to abuse terminology and refer to any $G$-set isomorphic to $G/H$ for some $H$ as an \emph{orbit}. This is important enough that we record it in the corollary below.

\begin{corollary}
	\label{G-sets are disjoint unions of coset spaces}
	Any $G$-set $X$ can be written as a disjoint union (coproduct)
	\[
		X \cong \coprod_i G/H_i
	\]
	for subgroups $H_i \leq G$.
\end{corollary}

	Orbits act like ``atoms'' -- they are the indecomposable building blocks of $G$-sets. In this way, they play a role similar to that of \emph{elements} in the category of sets.

\begin{proposition}
	\label{orbits to orbits}
	Let $f \colon X \to Y$ be a morphism of $G$-sets.
	\begin{enumerate}[(a)]
		\item The image of each orbit of $X$ under $f$ is an orbit of $Y$.
		\item The preimage of each orbit of $Y$ is a union of orbits of $X$.
		\item Let $Y = \coprod_i Y_i$, where each $Y_i$ is an orbit of $Y$. The morphism $f$ is isomorphic (in the arrow category of $\GSet$) to the coproduct over $i$ of the restrictions $f|_{f^{-1}(Y_i)} \colon f^{-1}(Y_i) \to Y_i$.
	\end{enumerate}
\end{proposition}

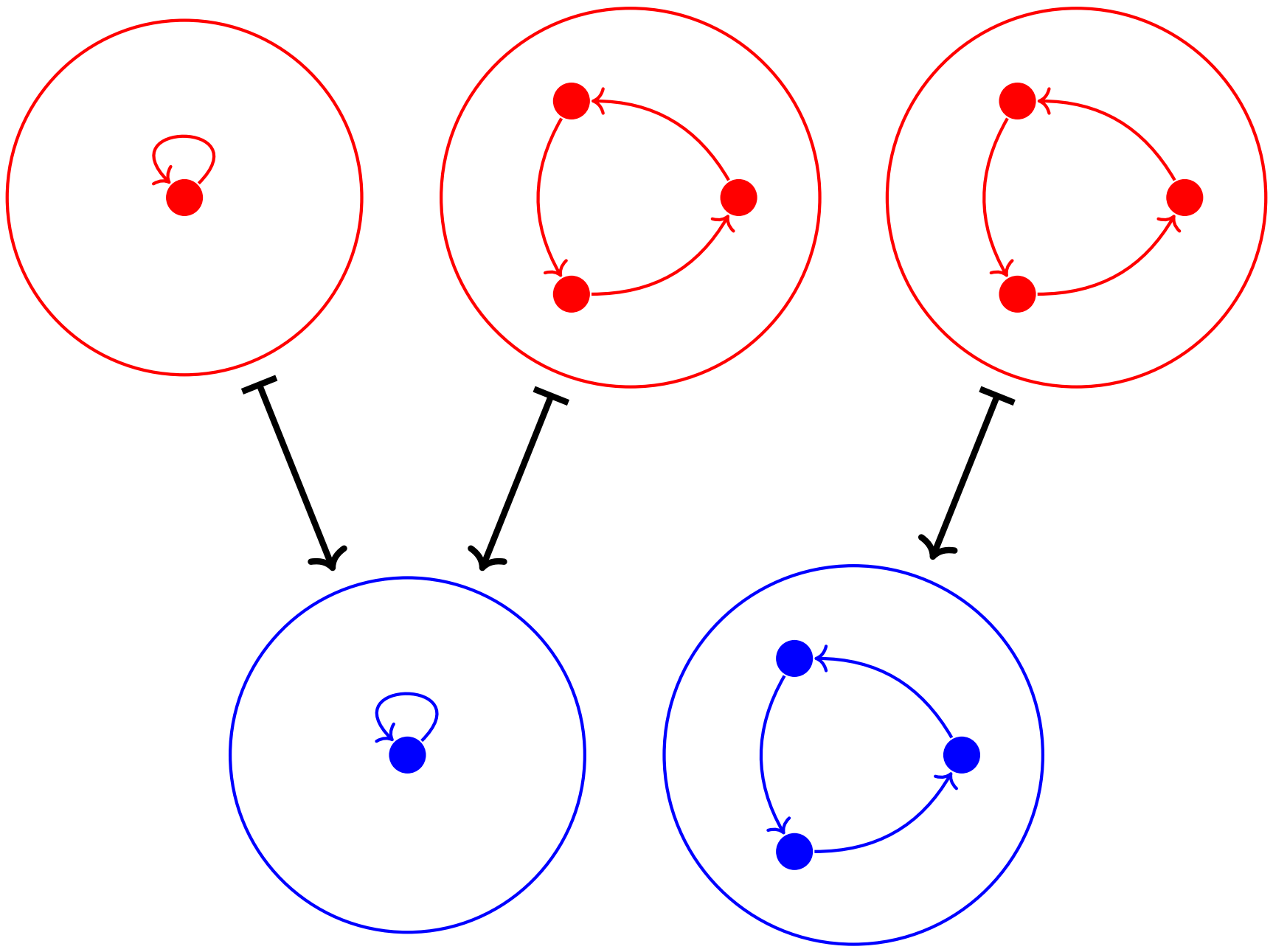
\begin{figure}[h]
	\begin{center}
		\begin{tikzpicture}
	\node[thick,blue,circle,draw] (B) at (2,0)
	{
		\begin{tikzpicture}
			\node[thick,fill,circle] (e) at (0:1) {};
			\node[thick,fill,circle] (g) at (120:1) {};
			\node[thick,fill,circle] (g2) at (240:1) {};
			\draw[->,thick] (e) to[bend right=30] (g);
			\draw[->,thick] (g) to[bend right=30] (g2);
			\draw[->,thick] (g2) to[bend right=30] (e);
			\draw[white] (180:1) circle (0.1);
		\end{tikzpicture}
	};
	\node[thick,blue,circle,draw] (A) at (-2,0)
	{
		\begin{tikzpicture}
			\node[thick,fill,circle] (e) at (0,0) {};
			\draw[white] (0,0) circle (1);
			\draw[->,thick] (e) to[loop] (e);
		\end{tikzpicture}
	};
	\node[thick,red,circle,draw] (Z) at (4,5)
	{
		\begin{tikzpicture}
			\node[thick,fill,circle] (e) at (0:1) {};
			\node[thick,fill,circle] (g) at (120:1) {};
			\node[thick,fill,circle] (g2) at (240:1) {};
			\draw[->,thick] (e) to[bend right=30] (g);
			\draw[->,thick] (g) to[bend right=30] (g2);
			\draw[->,thick] (g2) to[bend right=30] (e);
			\draw[white] (180:1) circle (0.1);
		\end{tikzpicture}
	};
	\node[thick,red,circle,draw] (Y) at (0,5)
	{
		\begin{tikzpicture}
			\node[thick,fill,circle] (e) at (0:1) {};
			\node[thick,fill,circle] (g) at (120:1) {};
			\node[thick,fill,circle] (g2) at (240:1) {};
			\draw[->,thick] (e) to[bend right=30] (g);
			\draw[->,thick] (g) to[bend right=30] (g2);
			\draw[->,thick] (g2) to[bend right=30] (e);
			\draw[white] (180:1) circle (0.1);
		\end{tikzpicture}
	};
	\node[thick,red,circle,draw] (X) at (-4,5)
	{
		\begin{tikzpicture}
			\node[thick,fill,circle] (e) at (0,0) {};
			\draw[white] (0,0) circle (1);
			\draw[->,thick] (e) to[loop] (e);
		\end{tikzpicture}
	};
	\draw[ultra thick,|->,shorten >=5,shorten <=5] (X) -- (A);
	\draw[ultra thick,|->,shorten >=5,shorten <=5] (Y) -- (A);
	\draw[ultra thick,|->,shorten >=5,shorten <=5] (Z) -- (B);
\end{tikzpicture}
	\end{center}
	\caption{A schematic illustration of a morphism of $C_3$-sets $f \colon X \to Y$. A more precise illustration of the function $f$ would label the elements of domain and codomain and indicate the image of each element of the domain. Instead, this illustration shows the ``essential structure'' of the morphism $f$, i.e. which orbits of the domain are mapped on which of the codomain. The orbits of $X$ are isomorphic to $C_3/C_3$, $C_3/e$, and $C_3/e$ (top row, in red), and the orbits of $Y$ are isomorphic to $C_3/C_3$ and $C_3/e$ (bottom row, in blue). The colored arrows indicate the action of a generator of $C_3$.}
	\label{morphism of C3 sets}
\end{figure}

See \Cref{morphism of C3 sets} for an illustration. In light of this proposition, understanding morphisms between arbitrary $G$-sets can be reduced to understanding morphisms between orbits, i.e., $G$-sets isomorphic to $G/H$.

Recall that for subgroups $K$ and $H$ of $G$, we say that $K$ is \emph{subconjugate} to $H$ if there exists $g \in G$ such that $gKg^{-1} \subseteq H$.

\begin{proposition}
	\label{subconjugate and quotients}
	Morphisms $G/K \to G/H$ are in bijection with cosets $gH \in G/H$ such that $g^{-1}Kg \subseteq H$ via the bijection that sends $gH$ to the morphism
	\[
		xK \mapsto xgH.
	\]
	In other words, a morphism $G/K \to G/H$ is uniquely determined by the image of the trivial coset $eK$, and this image must be a coset $gH$ such that $g^{-1}Kg \subseteq H$.
\end{proposition}

In particular, there exists a morphism $G/K \to G/H$ if and only if $K$ is subconjugate to $H$. In the case that $K$ is a subgroup of $H$, there is a canonical quotient morphism $q^H_K \colon G/K \to G/H$ corresponding to the trivial coset $eH$, i.e. $q^H_K(gK) = gH$.

More generally, we can understand morphisms in $\GSet$ out of a single orbit as fixed points.

\begin{definition}
	Let $X$ be a $G$-set and let $H$ be a subgroup of $G$. The \emph{$H$-fixed points of $X$} are the set
	\[
		X^H := \{x \in X \mid \forall h \in H, h \cdot x = x\}.
	\]
\end{definition}

Equivalently, an $H$-fixed point is a point $x \in X$ such that the stabilizer subgroup of $x$ contains $H$.

\begin{proposition}
	\label{fixed points as morphisms from cosets}
	Let $H$ be a subgroup of $G$, and let $X$ be a $G$-set. Morphisms of $G$-sets $G/H \to X$ are in bijection with elements of $X^H$:
	\[
		\GSet(G/H,X) \cong X^H.
	\]
\end{proposition}

We end with a structure theorem for morphisms in the category $\GSet$.

\begin{theorem}
	\label{fold after quotient}
	Any morphism $f \colon X \to Y$ of $G$-sets is isomorphic (in the arrow category) to a coproduct of quotients followed by a coproduct of folds.
\end{theorem}

\begin{proof}
	By \cref{orbits to orbits}, we may reduce to the case that $Y$ has a single orbit, i.e. up to isomorphism we have $Y = G/H$ for some subgroup $H$ of $G$. Then choosing an isomorphism $X \cong \coprod_i G/K_i$, we have reduced to the case of a morphism of the form
	\[
		(f_i)_i \colon \coprod_i G/K_i \to G/H,
	\]
	which is equal to the composite
	\[
		\coprod_i G/K_i \xrightarrow{\coprod_i f_i} \coprod_i G/H \xrightarrow{\nabla} G/H.
	\]
	Note that this is a fold following a coproduct of morphisms $f_i \colon G/K_i \to G/H$.

	Thus, we have reduced to the case that $f \colon G/K \to G/H$ is a map between single orbits. It remains to show that $f$ is isomorphic to a quotient. Let $g \in G$ be such that $f(eK) = gH$; in particular we have $g^{-1}Kg \subseteq H$ from \Cref{subconjugate and quotients}. Then let $c \colon G/g^{-1}Kg \to G/K$ be the morphism (by \Cref{subconjugate and quotients}) corresponding to the coset $g^{-1}K$, i.e. $c(xg^{-1}Kg) = xg^{-1}K$. $c$ is an isomorphism for cardinality reasons. The quotient $q^H_{g^{-1}Kg}$ and $c$ fit into a square
	\[
		\begin{tikzcd}
			G/K \ar[d, "f"'] & G/g^{-1}K g \ar[l, "c"'] \ar[d, "q^H_{g^{-1}Kg}"] \\
			G/H & G/H \ar[l, "\id"]
		\end{tikzcd}
	\]
	which one verifies commutes by computing the image of the trivial coset of $G/g^{-1}Kg$ along both composites. This witnesses the desired isomorphism in the arrow category between $f$ and $q^H_{g^{-1}Kg}$.
\end{proof}

\section{Lifting Properties in \texorpdfstring{$G$}{G}-Sets}
\label{section: lifting properties in G-sets}

Fix a finite group $G$. We apply the strategy outlined in \cref{section: the strategy} to the category of $G$-sets. Although we cannot fully execute the classification here, we will make substantially more progress in \cref{section: Cp-sets} when $G$ is a group of prime order.

By \cref{fold after quotient}, any morphism of $G$-sets is isomorphic to a coproduct of quotients followed by a coproduct of fold maps.

\begin{definition}
	\label{def quotient and fold}
	Let $K \leq H$ be subgroups of $G$.
	\begin{enumerate}[(a)]
		\item Let
		      \(
		      q^H_K \colon G/K \to G/H
		      \)
		      denote the quotient map from $G/K$ to $G/H$.
		\item For any set $I$, let
		      \(
		      \nabla_{H,I} \colon \coprod_{i \in I} G/H \to G/H
		      \)
		      denote the fold map: the morphism out of the coproduct induced by identity maps $G/H \to G/H$. When $H$ is clear from context, we will often write $\nabla_{I}$ instead of $\nabla_{H,I}$. We will also often abuse notation and write $\nabla_{H,n}$ to denote $\nabla_{H,I}$ for some implicit indexing set $I$ with $\lvert I \rvert = n$, understanding that $n$ may be an infinite cardinal if $I$ is infinite.
	\end{enumerate}
\end{definition}

Given that every morphism is built under coproduct and composition out of quotients and folds by \cref{fold after quotient}, a good starting point for our classification is to find the right-lifting classes and the vacuous right-lifting classes of a single quotient $q^H_K$ and a single fold $\nabla_{H,n}$. Then we may combine these using \cref{lifting class of coproduct,lifting class of union}.

We introduce in \cref{table:morphisms} abbreviations for classes of morphisms that will appear frequently.

\begin{table}[h!]
	\centering
	\begin{tblr}{colspec={lcl}, rowsep=2pt}
		\toprule
		Class                        & Category       & Description                                            \\ \midrule
		$\inj$                       & $\Set$         & Injections                                             \\
		$\surj$                      & $\Set$         & Surjections                                            \\
		$\bij$                       & $\Set$         & Bijections                                             \\
		$\init$                      & $\Set$         & Functions with domain $\emptyset$                      \\
		$\inj^H$                     & $\GSet$        & Morphisms $f$ such that $f^H \in \inj$                 \\
		$\surj^H$                    & $\GSet$        & Morphisms $f$ such that $f^H \in \surj$                \\
		$\init^H$                    & $\GSet$        & Morphisms $f$ such that $f^H \in \init$                \\
		$\pb_K^H$                    & $\GSet$        & \Cref{def:pb}                                          \\
		$\homog$                     & $\GSet$        & \Cref{def:homog}                                       \\
		$\all$                       & $\GSet$        & All morphisms                                          \\
		$\{\id_\emptyset\}^H$ & $\GSet$ & Morphisms $f$ such that $f^H = \id_{\emptyset}$ \\
		\bottomrule
	\end{tblr}
	\caption{Classes of morphisms in $\GSet$. Here, $K$ and $H$ are subgroups of $G$ with $K \leq H$. }
	\label{table:morphisms}
\end{table}

\begin{definition}
	\label{def:pb}
	Let $\pb_K^H$ be the class of all morphisms $f : X \to Y$ in $\GSet$ such that
	\[\begin{tikzcd}[ampersand replacement=\&]
			{X^H} \& {X^K} \\
			{Y^H} \& {Y^K}
			\arrow[hook, from=1-1, to=1-2]
			\arrow["{f^H}"', from=1-1, to=2-1]
			\arrow["{f^K}", from=1-2, to=2-2]
			\arrow[hook, from=2-1, to=2-2]
		\end{tikzcd}\]
	is a pullback diagram in $\Set$. This is equivalent to the statement
	\begin{equation}
		\label{equivalent characterization of pb}
		f^{-1}(Y^H) \cap X^K = X^H.
	\end{equation}
\end{definition}

It is also worth noting here that we always have $f(X^H) \subseteq Y^H$; an equivariant function sends $H$-fixed points of $X$ to $H$-fixed points of $Y$. We will often use this containment in the form $X^H \subseteq f^{-1}(Y^H)$. 

\begin{definition}
	\label{def:homog}
	A $G$-set $X$ is said to be \emph{homogeneous} if the orbits of $X$ are pairwise-isomorphic (equivalently, the stabilizers of elements of $X$ are pairwise-conjugate). A morphism $f : X \to Y$ in the category of $G$-sets is said to be \emph{homogeneous} if, for all orbits $O \subseteq Y$, $f^{-1}(O)$ is a homogeneous $G$-set. The class of homogeneous morphisms in $\GSet$ will be denoted $\homog$.
\end{definition}

\begin{proposition}\label{prop:nabla boxslash}
	Let $n$ be a (possibly infinite) cardinal. For any subgroup $H \leq G$, we have
	\[\rlp{\{\nabla_{H,n}\}} = \begin{cases} \surj^H &: n = 0 \\ \all &: n = 1 \\ \inj^H &: n > 1 \end{cases}\]
\end{proposition}

\begin{proof}
	Consider the adjunction
	\[\begin{tikzcd}[ampersand replacement=\&]
			\GSet \\
			\Set,
			\arrow[""{name=0, anchor=center, inner sep=0}, "{R}", shift left=2, from=1-1, to=2-1]
			\arrow[""{name=1, anchor=center, inner sep=0}, "{L}", shift left=2, from=2-1, to=1-1]
			\arrow["\dashv"{anchor=center}, draw=none, from=1, to=0]
		\end{tikzcd}\]
	where the left adjoint $L \colon \Set \to \GSet$ is given by $L(A) = A \times G/H$ (where we view a set $A$ as a $G$-set with trivial $G$-action). The right adjoint $R \colon \GSet \to \Set$ is given by $R(X) = X^H$. Then $\nabla_{\!H,n}$ is the image under the left adjoint of the unique function $\nabla_{\!n} \colon I \to \ast$, where $\lvert I \rvert = n$. Now
		\begin{align*}
			\rlp{\{\nabla_{H,n}\}}
			 & = \rlp{\{L(\nabla_{\!n})\}}
			\\
			 & = R^{-1}\left(\rlp{\{\nabla_{\!n}\}}\right)
			 & \textnormal{\Cref{proposition:adjunction boxslash}}
			\\
			 & =
			\begin{cases}
				R^{-1}(\surj) & \colon n = 0 \\
				R^{-1}(\all)  & \colon n = 1 \\
				R^{-1}(\inj)  & \colon n > 1 \\
			\end{cases}
			 & \textnormal{\cref{lifts of folds}}
			\\
			 & =
			\begin{cases}
				\surj^H & \colon n = 0  \\
				\all    & \colon n = 1  \\
				\inj^H  & \colon n > 1. \\
			\end{cases}
			\qedhere
		\end{align*}
\end{proof}

\begin{proposition}
	\label{lemma:vacuous lifting classes of nabla}
	Let $n$ be a (possibly infinite) cardinal. For any subgroup $H \leq G$, we have
		\[
			\{\nabla_{H,n}\}^\emptyset =
			\begin{cases}
				\{\id_\emptyset\}^H & \colon n = 0 \\
				\init^H             & \colon n > 0
			\end{cases}
		\]
\end{proposition}

\begin{proof}
	Consider first the case $n = 0$. A morphism $g \colon X \to Y$ lies in $\{\nabla_{H,0}\}^\emptyset$ if and only if there exist no commutative squares in $\GSet$ of the form 
	\[
		\begin{tikzcd}
			\emptyset \ar[r] \ar[d, "\nabla_{H,0}"'] & X \ar[d, "g"] \\
			G/H \ar[r] & Y.
		\end{tikzcd}
	\]
	Since $\emptyset$ is the initial object in $\GSet$, this is equivalent to there existing no morphism $G/H \to Y$. By \cref{fixed points as morphisms from cosets}, morphisms $G/H \to Y$ correspond to fixed points $Y^H$. So, $g \in \{\nabla_{\!H,0}\}^\emptyset$ if and only if $Y^H = \emptyset$. Then we simply note that $Y^H = \emptyset$ is equivalent to $g^H = \id_\emptyset$, i.e.\ $g \in \{\id_\emptyset\}^H$.

	When $n > 0$, let $I$ be an indexing set with $\lvert I \rvert = n$. A morphism $g \colon X \to Y$ lies in $\{\nabla_{H,n}\}^\emptyset$ if and only if there exist no commutative squares in $\GSet$ of the form 
	\[
		\begin{tikzcd}
			\coprod_{i \in I} G/H \ar[r] \ar[d, "\nabla_{H,n}"'] & X \ar[d, "g"] \\
			G/H \ar[r] & Y.
		\end{tikzcd}
	\]
	Using \cref{fixed points as morphisms from cosets} liberally, such commutative diagrams correspond bijectively with tuples $(x_i)_{i \in I}$ of elements of $X^H$ such that $g(x_i) = g(x_j)$ for all $i,j \in I$. This correspondence is given as follows:
	\begin{itemize}
		\item $x_i$ is the image of the identity coset $eH$ in the $i$th copy of $G/H$ in the domain of $\nabla_{H,n}$.
		\item $g(x_i)$ is the image of the identity coset $eH$ in the codomain of $\nabla_{H,n}$. The requirement $x_i \in X^H$ implies $g(x_i) \in Y^H$.
		\item Since each basepoint $eH$ in the domain of $\nabla_{H,n}$ is sent to the basepoint $eH$ in the codomain, the condition that the square commutes is equivalent to the condition that $g(x_i) = g(x_j)$ for all $i,j$.
		\end{itemize}
	Note that there are no such tuples if and only if $X^H = \emptyset$ (else, choose any $x \in X^H$ and use the constant tuple with $x_i = x$ for all $i$). This is equivalent to $g^H \in \init^H$, as desired.
\end{proof}

\begin{lemma}
	\label{fixed points to pullback is always injective}
	Let $f \colon X \to Y$ be a morphism in $\GSet$, and let $K \leq H$ be subgroups of $G$. Let $P$ be the pullback of the diagram 
	\[
		\begin{tikzcd}
			P
			\ar[dr, phantom, "\lrcorner" description, very near start]
			\ar[r, "p"]
			\ar[d]
			&
			X^K
			\ar[d, "f^K"]
			\\
			Y^H
			\ar[r, hook]
			&
			Y^K.
		\end{tikzcd}
	\]
	The canonical map $\phi \colon X^H \to P$ determined by $f^H \colon X^H \to Y^H$ and the inclusion $i \colon X^H \hookrightarrow X^K$ is injective. 
\end{lemma}
\begin{proof}
	We have $i = p \circ \phi$, and since $i$ is injective, $\phi$ must be injective as well.
\end{proof}

\begin{proposition}\label{proposition:qHKboxslash}
	For any subgroups $K \leq H \leq G$, we have
	\[
		\rlp{\{q^H_K\}} = \pb_K^H.
	\]
\end{proposition}

\begin{proof}
	Let $P$ be the following pullback in $\Set$:
	\[
		\begin{tikzcd}
			{P}
			\ar[dr, phantom, "\lrcorner" description, very near start]
			\ar[r]
			\ar[d]
			&
			{\GSet(G/K,X)}
			\ar[d, "f_*"]
			\\
			{\GSet(G/H,Y)}
			\ar[r, hook, "(q^H_K)^*"']
			&
			{\GSet(G/K,Y)}.
		\end{tikzcd}
	\]
	so that $P$ consists of pairs of maps $\alpha : G/K \to X$ and $\beta : G/H \to Y$ such that $f \circ \alpha = \beta \circ q^H_K$. Note that such pairs of maps correspond precisely to commutative squares
	\begin{equation}
		\label{rlp against quotient}
		\begin{tikzcd}
			G/K \ar[d, two heads, "q^H_K"']
			\ar[r, "\alpha"]
			&
			X
			\ar[d, "f"]
			\\
			G/H
			\ar[r, "\beta"']
			&
			Y.
		\end{tikzcd}
	\end{equation}

	Let $\phi \colon \GSet(G/H,X) \to P$ be the canonical morphism determined by $(q^H_K)^*$ and $f_*$, as in the diagram below:
	\[
		\begin{tikzcd}
			\GSet(G/H,X)
			\ar[drr, "(q^H_K)^*", bend left=20]
			\ar[ddr, "f_*"', bend right=30]
			\ar[dr, "\phi" description, dashed]
			\\
			&
			{P}
			\ar[dr, phantom, "\lrcorner" description, very near start]
			\ar[r]
			\ar[d]
			&
			{\GSet(G/K,X)}
			\ar[d, "f_*"]
			\\
			&
			{\GSet(G/H,Y)}
			\ar[r, "(q^H_K)^*"']
			&
			{\GSet(G/K,Y)}.
		\end{tikzcd}
	\]
	A square \eqref{rlp against quotient} admits a lift if and only if the corresponding element $(\alpha,\beta) \in P$ is in the image of $\phi$. Thus, we have $f \in \{q^H_K\}^{\boxslash}$ if and only if $\phi$ is surjective.

	By \cref{fixed points as morphisms from cosets}, this diagram is isomorphic to
	\[
		\begin{tikzcd}
			X^H
			\ar[drr, hook, bend left=20]
			\ar[ddr, "f^H"', bend right=30]
			\ar[dr, "\phi" description]
			\\
			&
			{P}
			\ar[dr, phantom, "\lrcorner" description, very near start]
			\ar[r]
			\ar[d]
			&
			{X^K}
			\ar[d, "f^K"]
			\\
			&
			{Y^H}
			\ar[r, hook]
			&
			{Y^K}.
		\end{tikzcd}
	\]
	Now by \cref{fixed points to pullback is always injective}, $\phi$ is always injective. Therefore, $f \in \{q^H_K\}^{\boxslash}$ if and only if $\phi$ is an isomorphism. In other words, $f \in \{q^H_K\}^{\boxslash}$ if and only if
	\[
		\begin{tikzcd}
			{X^H}
			\ar[r, hook]
			\ar[d, "f^H"']
			&
			{X^K}
			\ar[d, "f^K"]
			\\
			{Y^H}
			\ar[r, hook]
			&
			{Y^K}
		\end{tikzcd}
	\]
	is a pullback diagram in $\Set$. This is equivalent to $f \in \pb^H_K$ by definition.
\end{proof}

\begin{proposition}
	\label{vacuous lifting of quotients}
	For any subgroups $K \leq H \leq G$, we have
	\[
		\{q^H_K\}^\emptyset = \init^H \cap \pb^H_K.
	\]
\end{proposition}
\begin{proof}
	A morphism $g \colon X \to Y$ lies in $\{q^H_K\}^{\emptyset}$ if and only if there exist no commutative diagrams of the form 
	\[
		\begin{tikzcd}
			G/K \ar[r] \ar[d, "q^H_K"'] & X \ar[d, "g"] \\
			G/H \ar[r] & Y.
		\end{tikzcd}
	\]
	By \cref{fixed points as morphisms from cosets}, such commutative squares are in bijection with elements $x \in g^{-1}(Y^H) \cap X^K$. Therefore,
		\[g \in \{q^H_K\}^\emptyset \iff g^{-1}(Y^H) \cap X^K = \emptyset.\]
		Since we always have $X^H \subseteq g^{-1}(Y^H) \cap X^K$, this is equivalent to $X^H = \emptyset$ and $g^{-1}(Y^H) \cap X^K = X^H$. The former is the definition of $g \in \init^H$, and the latter is equivalent to $g \in \pb^H_K$ by \eqref{equivalent characterization of pb}.
\end{proof}

\section{Lifting Classes in Groups of Prime Order}
\label{section: Cp-sets}

To go any further, we must specialize to the case that $G$ has prime order. Fix a prime $p$ and write $G = C_p$. Importantly, this means that there are only two subgroups to consider: the trivial subgroup $e$ and the whole group $C_p$. In particular, every $C_p$-set is isomorphic to a coproduct of copies of $C_p/e$ and $C_p/C_p$ by \cref{G-sets are disjoint unions of coset spaces}. We call an orbit isomorphic to $C_p/e$ \emph{free}; an orbit isomorphic to $C_p/C_p$ is a fixed point.

A new class of morphisms which will play an important role in the following computations is $\{\id_\emptyset\}^{C_p}$, the class of all morphisms $f : X \to Y$ such that $X^{C_p} = \emptyset$ and $Y^{C_p} = \emptyset$. Equivalently, $\{\id_\emptyset\}^{C_p} = \{f : f^{C_p} = \id_\emptyset\}$.

To simplify notation in this context, we will write $\pb$ in place of $\pb_e^{C_p}$ for brevity. We can also give simpler characterizations of the classes $\pb$ and $\homog$ for groups of prime order. Recall that for any $f \colon X \to Y$ in $\CpSet$, we have $X^{C_p} \subseteq f^{-1}(Y^{C_p})$.
\begin{lemma}\
	\label{easier pb and homog for Cp}
	\begin{enumerate}[(a)]
		\item A morphism $f \colon X \to Y$ lies in $\pb$ if and only if $X^{C_p} = f^{-1}(Y^{C_p})$.
		\item A morphism $f \colon X \to Y$ lies in $\homog$ if and only if, for all $y \in Y^{C_p}$, $f^{-1}(\{y\})$ is either contained in or disjoint from $X^{C_p}$.
	\end{enumerate}
\end{lemma}

\subsection{Right-Lifting Classes of Morphisms with Transitive Codomain}

Let $f \colon X \to C_p/e$ be a morphism of $C_p$-sets. Writing $X \cong \coprod_{i \in I} X_i$ as a coproduct of its orbits $\{X_i \mid i \in I\}$, we see that $f$ is determined by a collection of morphisms $f_i \colon X_i \to C_p/e$, where $X_i$ is isomorphic to $C_p/H$ for a subgroup $H$. Since there are no morphisms $C_p/C_p \to C_p/e$ in $\CpSet$, each $X_i$ must be isomorphic to $C_p/e$, and we may choose the isomorphism such that $f_i \colon X_i \cong C_p/e \to C_p/e$ is an identity. In other words, up to precomposition with an automorphism of $X$, $f$ is a fold map $\nabla_{e,n}$ for some (possibly infinite) cardinal $n$.
The right-lifting classes and vacuous right-lifting classes with respect to these fold maps $\nabla_{e,n}$ are displayed in \cref{table:rlp for folds}.

\begin{table}[h]
	\centering
	\begin{tblr}{colspec={rll}, cells={mode=math}, colsep=6pt}
		\toprule
		n  & \{\nabla_{e,n}\}^\emptyset & \rlp{\{\nabla_{e,n}\}} \\ \midrule
		0  & \{\id_\emptyset\}          & \surj^e                \\ 
		1  & \init^e                    & \all                   \\ 
		>1 & \init^e                    & \inj^e                 \\
		\bottomrule
	\end{tblr}
	\vspace{1em}
	\caption{The right-lifting classes with respect to folds of free orbits in $\CpSet$. The right column comes from \cref{prop:nabla boxslash}, and the left column from \cref{lemma:vacuous lifting classes of nabla}. }
	\label{table:rlp for folds}
\end{table}

Next, let $f : X \to C_p/C_p$ be a morphism of $C_p$-sets. Then, up to isomorphism, we have $X = \coprod_{i \in I} C_p/e \amalg \coprod_{j \in J} C_p/C_p$, and $f$ is uniquely determined since $C_p/C_p$ is the terminal object. One can check that the right-lifting classes and vacuous right-lifting classes are as given in \cref{table:rlp for terminal codomain}. In the next two propositions, we prove that the rows of \cref{table:rlp for terminal codomain} not already covered by results in \cref{section: lifting properties in G-sets} are correct.

\begin{table}[h]
	\centering
	\begin{tblr}{colspec={rrlll}, cells={mode=math}, colsep=5pt}
		\toprule
		n      & m      & \{f\}^\emptyset         & \rlp{\{f\}}            & \text{Justification}                                                     \\
		\midrule
		0      & 0      & \{\id_\emptyset\}^{C_p} & \surj^{C_p}            & \text{\cref{prop:nabla boxslash,lemma:vacuous lifting classes of nabla}} \\
		0      & 1      & \init^{C_p}             & \all                   & \text{\cref{prop:nabla boxslash,lemma:vacuous lifting classes of nabla}} \\
		0      & >1     & \init^{C_p}             & \inj^{C_p}             & \text{\cref{prop:nabla boxslash,lemma:vacuous lifting classes of nabla}} \\
		1      & 0      & \init^{C_p} \cap \pb    & \pb                    & \text{\cref{vacuous lifting of quotients,proposition:qHKboxslash}}       \\
		>1     & 0      & \init^{C_p} \cap \pb    & \pb \cap \inj^{C_p}    & \text{\cref{n>1 m=0 in table 3}}                                         \\
		\geq 1 & \geq 1 & \init^{C_p}             & \homog \cap \inj^{C_p} & \text{\cref{n >= 1 m >= 1 in table 3}}                                   \\
		\bottomrule
	\end{tblr}
	\vspace{1em}
	\caption{The right-lifting classes and the vacuous right-lifting classes of morphisms $X \to C_p/C_p$ in $\CpSet$. We decompose $X$ as
		\[
			X \cong \coprod_{i \in I} C_p/e \amalg \coprod_{j \in J} C_p/C_p
		\]
		with $|I| = n$ and $|J| = m$,
		and note that the morphism is uniquely determined because its codomain is terminal.}
	\label{table:rlp for terminal codomain}
\end{table}

\begin{proposition}
	\label{n>1 m=0 in table 3}
	Let $I$ be an indexing set with $|I| = n$, where $n$ is some (possibly infinite) cardinal. Assume $n > 1$. Let $f$ be the unique morphism to the terminal object in $\CpSet$ with domain $\coprod_{i \in I} C_p/e$:
	\[
		f \colon \coprod_{i \in I} C_p/e \to C_p/C_p.
	\]
	Then $\rlp{\{f\}} = \pb \cap \inj^{C_p}$ and $\{f\}^\emptyset = \init^{C_p} \cap \pb$.
\end{proposition}

\begin{proof}
	As in \cref{lemma:vacuous lifting classes of nabla}, commutative squares
	\[
		\begin{tikzcd}
			\coprod_{i \in I} C_p/e \ar[d, "f"'] \ar[r] & X \ar[d, "g"] \\
			C_p / C_p \ar[r] & Y
		\end{tikzcd}
	\]
	correspond bijectively with tuples $(x_i)_{i \in I}$ of elements in $g^{-1}(Y^{C_p})$ such that $g(x_i) = g(x_j)$ for all $i, j \in I$; the element $x_i$ corresponds to the image under the top horizontal morphism of the trivial coset $e \in C_p/e$ in the $i$th copy of $C_p/e$.

	We first consider $\{f\}^\emptyset$. There are no such squares if and only if there are no tuples $(x_i)_{i \in I}$ as in the previous paragraph, and there are no such tuples if and only if $g^{-1}(Y^{C_p}) = \emptyset$.
		Since $X^{C_p} \subseteq g^{-1}(Y^{C_p})$, this is equivalent to having $g^{-1}(Y^{C_p}) = X^{C_p}$ and $X^{C_p} = \emptyset$. Thus, by \cref{easier pb and homog for Cp}, $\{f\}^\emptyset = \init^{C_p} \cap \pb$.

	We next consider $\rlp{\{f\}}$. By \cref{fixed points as morphisms from cosets}, a lift $C_p/C_p \to X$ corresponds to a choice of $x \in X^{C_p}$ such that $x_i = x$ for all $i \in I$ (in particular, a lift is unique if it exists). To have $g \notin \rlp{\{f\}}$, there must exist an element $y \in Y^{C_p}$ and a tuple $(x_i)_{i \in I}$ of elements in $g^{-1}(y)$ such that either:
		\begin{enumerate}[(i)]
			\item $x_i \notin X^{C_p}$ for some $i \in I$, or
			\item $x_i \neq x_j$ for some $i, j \in I$.
		\end{enumerate}
		In other words, we have $g \in \rlp{\{f\}}$ if and only if, for all $y \in Y^{C_p}$, the set $g^{-1}(y)$ is either empty or a singleton subset of $X^{C_p}$. This is equivalent to saying that $g^{-1}(Y^{C_p}) \subseteq X^{C_p}$, i.e. $g \in \pb$, and that $g^{C_p} \colon X^{C_p} \to Y^{C_p}$ is injective, i.e. $g \in \inj^{C_p}$. Therefore, $\rlp{\{f\}} = \pb \cap \inj^{C_p}$, as desired.
\end{proof}

\begin{proposition}
	\label{n >= 1 m >= 1 in table 3}
	Let $I$ and $J$ be indexing sets with $|I| = n$ and $|J|=m$, where $n$ and $m$ are some (possibly infinite) cardinals. Assume $n \geq 1$ and $m \geq 1$. Let $f$ be the unique morphism to the terminal object in $\CpSet$ as follows:
	\[
		f \colon \coprod_{i \in I} C_p/e \amalg \coprod_{j \in J} C_p/C_p  \to C_p/C_p.
	\]
	Then $\rlp{\{f\}} = \homog \cap \inj^{C_p}$ and $\{f\}^\emptyset = \init^{C_p}$.
\end{proposition}

\begin{proof}
	Similar to previous arguments, we see that commutative squares
	\[
		\begin{tikzcd}
			\coprod_{i \in I} C_p/e \amalg \coprod_{j \in J} C_p/C_p  \ar[d, "f"']
			\ar[r] & X \ar[d, "g"] \\
			C_p/C_p \ar[r] & Y
		\end{tikzcd}
	\]
	correspond bijectively with pairs of tuples $(x_i)_{i \in I}$ and $(x_j')_{j \in J}$ such that
		\begin{enumerate}[(i)]
			\item $x_j' \in X^{C_p}$ for all $j \in J$,
			\item all the elements have the same image under $g$: for all $i \in I$ and all $j \in J$, $g(x_i) = g(x'_j)$.
		\end{enumerate}
		Note that it follows that $x_i \in g^{-1}(Y^{C_p})$ for all $i \in I$.

	We first consider $\{f\}^\emptyset$. Note that if $X^{C_p} = \emptyset$, then there are no such squares, because there are no choices for the elements $x_j'$. Conversely, if $X^{C_p} \neq \emptyset$, then we can pick any element $x \in X^{C_p}$ and set $x_i = x_j' = x$ for all $i \in I$ and all $j \in J$, which gives a commutative square. Hence, $\{f\}^\emptyset = \init^{C_p}$.

		Next, consider $\rlp{\{f\}}$. A lift of such a square exists if and only if all of the elements are equal, i.e.\ for all $i \in I$ and all $j \in J$, $x_i = x'_j$ (lifts are unique when they exist). In other words, we have $g \in \rlp{\{f\}}$ if and only if, for all $y \in Y^{C_p}$, the fiber $g^{-1}(y)$ is either:
		\begin{enumerate}[(i)]
			\item disjoint from $X^{C_p}$ (in which case there are no choices for $x_j'$ such that $g(x_j') = y$), or
			\item a singleton subset of $X^{C_p}$ (in which case there is a unique choice for each $x_i$ and each $x_j'$).
		\end{enumerate}
		By \cref{easier pb and homog for Cp}, we obtain $\rlp{\{f\}} = \homog \cap \inj^{C_p}$, as desired.
\end{proof}

\subsection{Right-Lifting Classes of Morphisms with Arbitrary Codomain}

Beginning in this section and for the rest of the paper, the reader may find it helpful to follow along in the accompanying code \cite{code}.

Every morphism in $\CpSet$ is a coproduct of morphisms to transitive $C_p$-sets, so \cref{lifting class of coproduct} allows us to compute the possible values of $\rlp{\{f\}}$ for all morphisms $f$. Explicitly, $f = \coprod_i f_i$ with each $f_i$ a morphism with transitive codomain (i.e. a morphism whose lifting classes appear in \cref{table:rlp for folds,table:rlp for terminal codomain}), and
\begin{equation}
	\label{equation:rlp of coproduct}
	\rlp{\{f\}} = \rlp{\left\{\coprod_{i \in I} f_i \right\}}
	= \left(\bigcap_{i \in I} \{f_i\}^\boxslash\right) \cup \left(\bigcup_{i \in I} \{f_i\}^\emptyset\right).
\end{equation}

We begin by computing all possible unions of the classes $\{f_i\}^\emptyset$ as $f_i$ ranges over morphisms whose codomains are transitive. These vacuous lifting classes are drawn from those in \cref{table:rlp for folds,table:rlp for terminal codomain}:
\begin{equation}
	\label{equation:vacuous classes}
	\{\id_\emptyset\},\quad
	\{\id_\emptyset\}^{C_p},\quad
	\init^e,\quad
	\init^{C_p},\quad\text{ and  }\quad
	\init^{C_p} \cap \pb.
\end{equation}

\begin{figure}[h]
	\begin{tikzpicture}[xscale=2,yscale=1.2]
\node (v0) at (0, -2) {$ \{\id_\emptyset\} $};
\node (v1) at (-0.75, -1) {$ \init^e $};
\node (v2) at (0.75, -1) {$ \{\id_\emptyset\}^{C_p} $};
\node (v3) at (0, 0) {$ \init^e \cup \{\id_\emptyset\}^{C_p} $};
\node (v4) at (0, 1) {$ \init^{C_p} \cap \pb $};
\node (v5) at (0, 2) {$ \init^{C_p} $};
\draw[white, ultra thick, double=black, arrows={-Triangle[black,scale=0.25]}] (v0) -- (v1);
\draw[white, ultra thick, double=black, arrows={-Triangle[black,scale=0.25]}] (v0) -- (v2);
\draw[white, ultra thick, double=black, arrows={-Triangle[black,scale=0.25]}] (v1) -- (v3);
\draw[white, ultra thick, double=black, arrows={-Triangle[black,scale=0.25]}] (v2) -- (v3);
\draw[white, ultra thick, double=black, arrows={-Triangle[black,scale=0.25]}] (v3) -- (v4);
\draw[white, ultra thick, double=black, arrows={-Triangle[black,scale=0.25]}] (v4) -- (v5);
\node[fill=white] at (v0) {$ \{\id_\emptyset\} $};
\node[fill=white] at (v1) {$ \init^e $};
\node[fill=white] at (v2) {$ \{\id_\emptyset\}^{C_p} $};
\node[fill=white] at (v3) {$ \init^e \cup \{\id_\emptyset\}^{C_p} $};
\node[fill=white] at (v4) {$ \init^{C_p} \cap \pb $};
\node[fill=white] at (v5) {$ \init^{C_p} $};
\end{tikzpicture}
	\caption{Unions $\bigcup_i \{f_i\}^\emptyset$ of vacuous right-lifting classes of morphisms $f_i$ with transitive codomain. The classes in this poset are called \textit{collateral classes}.} 
	\label{fig:Cp empty classes}
\end{figure}
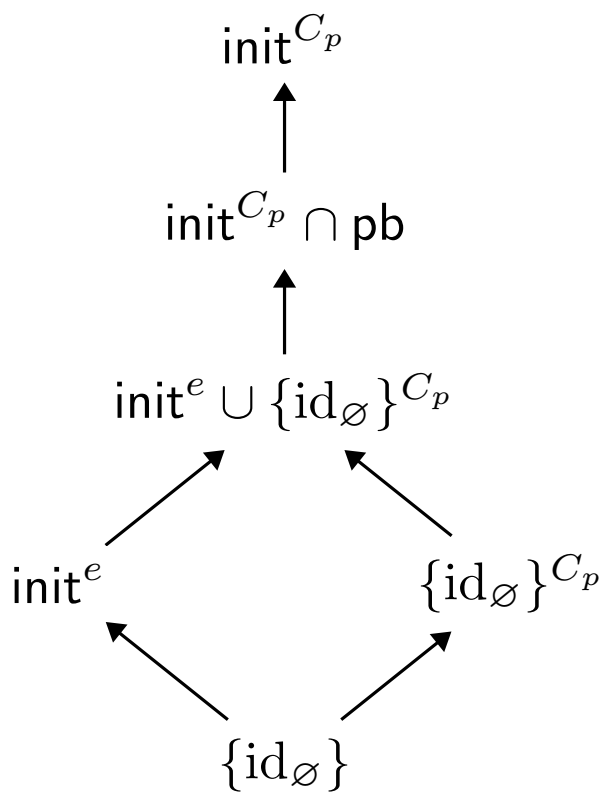

We call these five classes the \emph{primordial vacuous classes}. 
These classes are almost already closed under union; \Cref{fig:Cp empty classes} shows the unique union $\init^e \cup \{\id_\emptyset\}^{C_p}$ which does not already appear among the five vacuous lifting classes above.
Note that the inclusion $C_p/e \to C_p/e \amalg C_p/C_p$ lies in $\init^{C_p} \cap \pb$ but not in $\init^e \cup \{\id_\emptyset\}^{C_p}$, demonstrating that this new class is distinct from those in \eqref{equation:vacuous classes}.

Because it will be helpful below, we give the classes in \Cref{fig:Cp empty classes} a name.

\begin{definition}
	A \textit{collateral class} is one of the six classes $\{\id_\emptyset\}$, $\{\id_\emptyset\}^{C_p}$, $\init^e$, $\init^{C_p}$, $\init^e \cup \{\id_\emptyset\}^{C_p}$, and $\init^{C_p} \cap \pb$.
\end{definition}

In what follows, we will continue to introduce unusual names for certain classes of morphisms to connect them with the variable names in our code. We provide a ``cheat sheet'' in \Cref{fig:cheat sheet} which the reader may refer to for a quick reference.

\begin{strategy}
	\label{strategy}
	We would like to know that this picture is accurate, i.e.\ that \Cref{fig:Cp empty classes} is in fact the Hasse diagram of the poset of all possible unions of vacuous lifting classes. We have already argued abstractly that every possible union of vacuous lifting classes is one of the classes shown above, and that every arrow in the diagram depicts containment between classes. What remains to be shown is that there are no containments between classes not depicted in the diagram: for every non-containment $\mathcal P \nsubseteq \mathcal Q$ claimed by the diagram, we must produce a witness $f \in \mathcal P \setminus \mathcal Q$.

	This is a procedure we will have to execute frequently in what follows, to guarantee that the posets produced genuinely represent the containments between classes of morphisms in $\CpSet$. As a general strategy, we produce in \cref{appendix:proper containments} a collection of example morphisms and two tables (\cref{table:primordial proper containments,table:primordial vacuous proper containments}) listing which classes from \cref{table:morphisms} these examples lie in. Then we encode these tables in our code \cite{code} and computer verify that each claimed non-containment in a given poset has a witness. 
\end{strategy}

\bigskip

We now turn to computing all possible intersections of classes $\{f_i\}^\boxslash$ as $f_i$ ranges over morphisms with transitive codomain; this is the other piece of \cref{equation:rlp of coproduct} we must compute to determine all lifting classes of morphisms. The possible classes $\{f_i\}^\boxslash$ are drawn from \cref{table:rlp for folds,table:rlp for terminal codomain}:
\begin{equation}
	\label{equation:right lifting classes}
	\all,\
	\inj^e,\
	\surj^e,\
	\inj^{C_p},\
	\surj^{C_p},\
	\pb, \
	\pb \cap \inj^{C_p}, \text{ and }
	\homog \cap \inj^{C_p}.
\end{equation}
To compute intersections between these classes, we first determine the pairwise containments between the atomic pieces $\inj^e, \surj^e, \inj^{C_p}, \surj^{C_p}, \pb$, and $\homog$. We give these classes a name.

\begin{definition}
	A \emph{primordial class} is one of the six classes $\inj^e$, $\surj^e$, $\inj^{C_p}$, $\surj^{C_p}$, $\pb$, or $\homog$. A \textit{primordial intersection} is an intersection of some number of these classes.
\end{definition}

The containments between primordial classes are recorded in \cref{table:6x6containments}.

\begin{table}[h]
	\centering
	\begin{tblr}{colspec={c|*{6}{c}}, cells={mode=math}, colsep=5pt}
		\toprule
		\downarrow \subseteq \rightarrow & \inj^e & \inj^{C_p} & \surj^e & \surj^{C_p} & \pb    & \homog \\
		\midrule
		\inj^e                           & \cmark & \cmark     & \xmark  & \xmark      & \cmark & \cmark \\
		\inj^{C_p}                       & \xmark & \cmark     & \xmark  & \xmark      & \xmark & \xmark \\
		\surj^e                          & \xmark & \xmark     & \cmark  & \xmark      & \xmark & \xmark \\
		\surj^{C_p}                      & \xmark & \xmark     & \xmark  & \cmark      & \xmark & \xmark \\
		\pb                              & \xmark & \xmark     & \xmark  & \xmark      & \cmark & \cmark \\
		\homog                           & \xmark & \xmark     & \xmark  & \xmark      & \xmark & \cmark \\
		\bottomrule
	\end{tblr}
	\caption{Pairwise containments between the primordial classes $\inj^e$, $\surj^e$, $\inj^{C_p}$, $\surj^{C_p}$, $\pb$, and $\homog$.}
	\label{table:6x6containments}
\end{table}

Next, we note that $\inj^e \cap \surj^e = \bij^e$ is the class of isomorphisms, which is contained in all six primordial classes. There are also two nontrivial containments of intersections of primordial classes, described in the next two lemmas.

\begin{lemma}
	\label{primordial_relation_1}
	$\surj^e \cap \pb \subseteq \surj^{C_p}$.
\end{lemma}
\begin{proof}
	Suppose $f : X \to Y$ is a morphism which lies in both $\surj^e$ and $\pb$. Let $y \in Y^{C_p}$ be arbitrary. Since $f \in \surj^e$, there is some $x \in X$ such that $f(x) = y$. If $x \notin X^{C_p}$, then the orbit $O$ of $x$ is free and $f(O) = \{y\}$, contradicting the assumption that $f \in \pb$ (by \cref{easier pb and homog for Cp}). Thus, $x \in X^{C_p}$. Since $y$ was arbitrary, we have $f \in \surj^{C_p}$.
\end{proof}

\begin{lemma}
	\label{primordial_relation_2}
	$\surj^{C_p} \cap \homog \subseteq \pb.$
\end{lemma}

\begin{proof}
	Recall the descriptions of $\pb$ and $\homog$ in $\CpSet$ from \cref{easier pb and homog for Cp}. Suppose $f : X \to Y$ is a morphism which lies in both $\surj^{C_p}$ and $\homog$. Suppose for contradiction that there exists a free orbit $O \subseteq X$ such that $\{y\} := f(O) \subseteq Y^{C_p}$. Since $f \in \surj^{C_p}$, there is some $x \in X^{C_p}$ such that $f(x) = y$. Now $O \cup \{x\} \subseteq f^{-1}(\{y\})$, contradicting the assumption $f \in \homog$. Thus, $f^{-1}(Y^{C_p}) \subseteq X^{C_p}$. Since we always have $f^{-1}(Y^{C_p}) \supseteq X^{C_p}$, we see that these sets are equal, whence $f \in \pb$ by \cref{easier pb and homog for Cp}.
\end{proof}

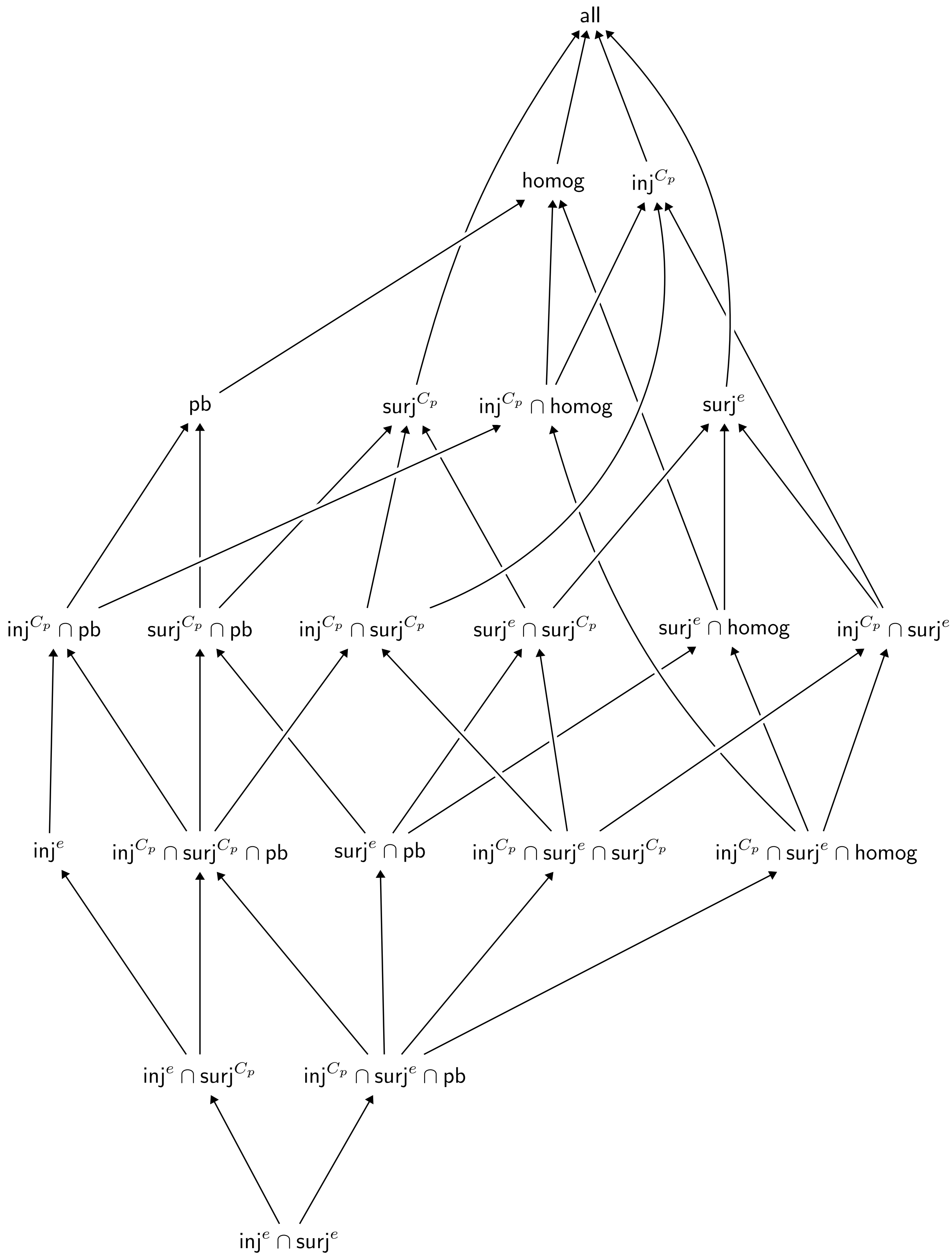
\begin{figure}[p]
	\begin{tikzpicture}[scale=3.5]
    \node (v0) at (0.5, 2.75) {$ \all $};
    \node (v1) at (0.787, 2) {$ \inj^{C_p} $};
    \node (v2) at (1.101, 1) {$ \surj^e $};
    \node (v3) at (-0.303, 1) {$ \surj^{C_p} $};
    \node (v4) at (0.334, 2) {$ \homog $};
    \node (v5) at (1.86, 0) {$ \inj^{C_p} \cap \surj^e $};
    \node (v6) at (-0.52, 0) {$ \inj^{C_p} \cap \surj^{C_p} $};
    \node (v7) at (0.3, 1) {$ \inj^{C_p} \cap \homog $};
    \node (v8) at (0.257, 0) {$ \surj^e \cap \surj^{C_p} $};
    \node (v9) at (1.101, 0) {$ \surj^e \cap \homog $};
    \node (v10) at (-1.25, 1) {$ \pb $};
    \node (v11) at (0.41, -1) {$ \inj^{C_p} \cap \surj^e \cap \surj^{C_p} $};
    \node (v12) at (1.513, -1) {$ \inj^{C_p} \cap \surj^e \cap \homog $};
    \node (v13) at (-1.905, 0) {$ \inj^{C_p} \cap \pb $};
    \node (v14) at (-1.25, 0) {$ \surj^{C_p} \cap \pb $};
    \node (v15) at (-1.926, -1) {$ \inj^e $};
    \node (v16) at (-1.25, -1) {$ \inj^{C_p} \cap \surj^{C_p} \cap \pb $};
    \node (v17) at (-0.443, -1) {$ \surj^e \cap \pb $};
    \node (v18) at (-1.25, -2) {$ \inj^e \cap \surj^{C_p} $};
    \node (v19) at (-0.422, -2) {$ \inj^{C_p} \cap \surj^e \cap \pb $};
    \node (v20) at (-0.85, -2.75) {$ \inj^e \cap \surj^e $};
    \draw[white, ultra thick, double=black, arrows={-Triangle[black,scale=0.25]}] (v20) -- (v19);
    \draw[white, ultra thick, double=black, arrows={-Triangle[black,scale=0.25]}] (v20) -- (v18);
    \draw[white, ultra thick, double=black, arrows={-Triangle[black,scale=0.25]}] (v19) -- (v17);
    \draw[white, ultra thick, double=black, arrows={-Triangle[black,scale=0.25]}] (v19) -- (v12);
    \draw[white, ultra thick, double=black, arrows={-Triangle[black,scale=0.25]}] (v19) -- (v11);
    \draw[white, ultra thick, double=black, arrows={-Triangle[black,scale=0.25]}] (v19) -- (v16);
    \draw[white, ultra thick, double=black, arrows={-Triangle[black,scale=0.25]}] (v17) -- (v9);
    \draw[white, ultra thick, double=black, arrows={-Triangle[black,scale=0.25]}] (v17) -- (v8);
    \draw[white, ultra thick, double=black, arrows={-Triangle[black,scale=0.25]}] (v17) -- (v14);
    \draw[white, ultra thick, double=black, arrows={-Triangle[black,scale=0.25]}] (v12) -- (v9);
    \draw[white, ultra thick, double=black, arrows={-Triangle[black,scale=0.25]}] (v12) -- (v5);
    \draw[white, ultra thick, double=black, arrows={-Triangle[black,scale=0.25]}] (v12) to[bend left=15] (v7);
    \draw[white, ultra thick, double=black, arrows={-Triangle[black,scale=0.25]}] (v9) -- (v2);
    \draw[white, ultra thick, double=black, arrows={-Triangle[black,scale=0.25]}] (v9) -- (v4);
    \draw[white, ultra thick, double=black, arrows={-Triangle[black,scale=0.25]}] (v11) -- (v8);
    \draw[white, ultra thick, double=black, arrows={-Triangle[black,scale=0.25]}] (v11) -- (v5);
    \draw[white, ultra thick, double=black, arrows={-Triangle[black,scale=0.25]}] (v11) -- (v6);
    \draw[white, ultra thick, double=black, arrows={-Triangle[black,scale=0.25]}] (v8) -- (v2);
    \draw[white, ultra thick, double=black, arrows={-Triangle[black,scale=0.25]}] (v8) -- (v3);
    \draw[white, ultra thick, double=black, arrows={-Triangle[black,scale=0.25]}] (v5) -- (v2);
    \draw[white, ultra thick, double=black, arrows={-Triangle[black,scale=0.25]}] (v5) -- (v1);
    \draw[white, ultra thick, double=black, arrows={-Triangle[black,scale=0.25]}] (v2) to[bend right=25] (v0);
    \draw[white, ultra thick, double=black, arrows={-Triangle[black,scale=0.25]}] (v18) -- (v16);
    \draw[white, ultra thick, double=black, arrows={-Triangle[black,scale=0.25]}] (v18) -- (v15);
    \draw[white, ultra thick, double=black, arrows={-Triangle[black,scale=0.25]}] (v16) -- (v14);
    \draw[white, ultra thick, double=black, arrows={-Triangle[black,scale=0.25]}] (v16) -- (v6);
    \draw[white, ultra thick, double=black, arrows={-Triangle[black,scale=0.25]}] (v16) -- (v13);
    \draw[white, ultra thick, double=black, arrows={-Triangle[black,scale=0.25]}] (v14) -- (v3);
    \draw[white, ultra thick, double=black, arrows={-Triangle[black,scale=0.25]}] (v14) -- (v10);
    \draw[white, ultra thick, double=black, arrows={-Triangle[black,scale=0.25]}] (v6) -- (v3);
    \draw[white, ultra thick, double=black, arrows={-Triangle[black,scale=0.25]}] (v6) to[bend right=40] (v1);
    \draw[white, ultra thick, double=black, arrows={-Triangle[black,scale=0.25]}] (v3) to[bend left=10] (v0);
    \draw[white, ultra thick, double=black, arrows={-Triangle[black,scale=0.25]}] (v15) -- (v13);
    \draw[white, ultra thick, double=black, arrows={-Triangle[black,scale=0.25]}] (v13) -- (v10);
    \draw[white, ultra thick, double=black, arrows={-Triangle[black,scale=0.25]}] (v13) -- (v7);
    \draw[white, ultra thick, double=black, arrows={-Triangle[black,scale=0.25]}] (v10) -- (v4);
    \draw[white, ultra thick, double=black, arrows={-Triangle[black,scale=0.25]}] (v7) -- (v4);
    \draw[white, ultra thick, double=black, arrows={-Triangle[black,scale=0.25]}] (v7) -- (v1);
    \draw[white, ultra thick, double=black, arrows={-Triangle[black,scale=0.25]}] (v4) -- (v0);
    \draw[white, ultra thick, double=black, arrows={-Triangle[black,scale=0.25]}] (v1) -- (v0);
    \node[fill=white] at (v0) {$ \all $};
    \node[fill=white] at (v1) {$ \inj^{C_p} $};
    \node[fill=white] at (v2) {$ \surj^e $};
    \node[fill=white] at (v3) {$ \surj^{C_p} $};
    \node[fill=white] at (v4) {$ \homog $};
    \node[fill=white] at (v5) {$ \inj^{C_p} \cap \surj^e $};
    \node[fill=white] at (v6) {$ \inj^{C_p} \cap \surj^{C_p} $};
    \node[fill=white] at (v7) {$ \inj^{C_p} \cap \homog $};
    \node[fill=white] at (v8) {$ \surj^e \cap \surj^{C_p} $};
    \node[fill=white] at (v9) {$ \surj^e \cap \homog $};
    \node[fill=white] at (v10) {$ \pb $};
    \node[fill=white] at (v11) {$ \inj^{C_p} \cap \surj^e \cap \surj^{C_p} $};
    \node[fill=white] at (v12) {$ \inj^{C_p} \cap \surj^e \cap \homog $};
    \node[fill=white] at (v13) {$ \inj^{C_p} \cap \pb $};
    \node[fill=white] at (v14) {$ \surj^{C_p} \cap \pb $};
    \node[fill=white] at (v15) {$ \inj^e $};
    \node[fill=white] at (v16) {$ \inj^{C_p} \cap \surj^{C_p} \cap \pb $};
    \node[fill=white] at (v17) {$ \surj^e \cap \pb $};
    \node[fill=white] at (v18) {$ \inj^e \cap \surj^{C_p} $};
    \node[fill=white] at (v19) {$ \inj^{C_p} \cap \surj^e \cap \pb $};
    \node[fill=white] at (v20) {$ \inj^e \cap \surj^e $};
\end{tikzpicture}
	\caption{The meet-semilattice generated by the primordial classes $\inj^e$, $\inj^{C_p}$, $\surj^e$, $\surj^{C_p}$, $\pb$, and $\homog$, from \cref{table:6x6containments} and the relations obtained in \cref{primordial_relation_1,primordial_relation_2}. The classes in this figure are called \textit{primordial intersections}. This figure is generated by computer, and differs from the output of our code \textit{only} by changing the positions of nodes for readability.}
	\label{fig:primordial_meet_semilattice}
\end{figure}

Now we can form the free meet-semilattice generated by the primordial classes $\inj^e$, $\surj^e$, $\inj^{C_p}$, $\surj^{C_p}$, $\pb$, and $\homog$, modulo the relations described above. All possible intersections of the primordial classes (the \textit{primordial intersections}) are shown in \Cref{fig:primordial_meet_semilattice}, which was obtained from our code \cite{code}. Using \cref{strategy}, we verify by computer that every non-containment claimed by \Cref{fig:primordial_meet_semilattice} is witnessed by an example morphism. In particular, the $21$ classes shown are distinct, and every containment shown is proper.

This semilattice is larger than the semilattice of all possible intersections $\bigcap_{i \in I} \rlp{\{f_i\}}$, because the class $\homog$ does not appear as the right-lifting class $\rlp{\{f_i\}}$ of any single morphism $f_i$ with transitive codomain \cref{equation:right lifting classes}, and it is not an intersection of classes from \cref{equation:right lifting classes}. 

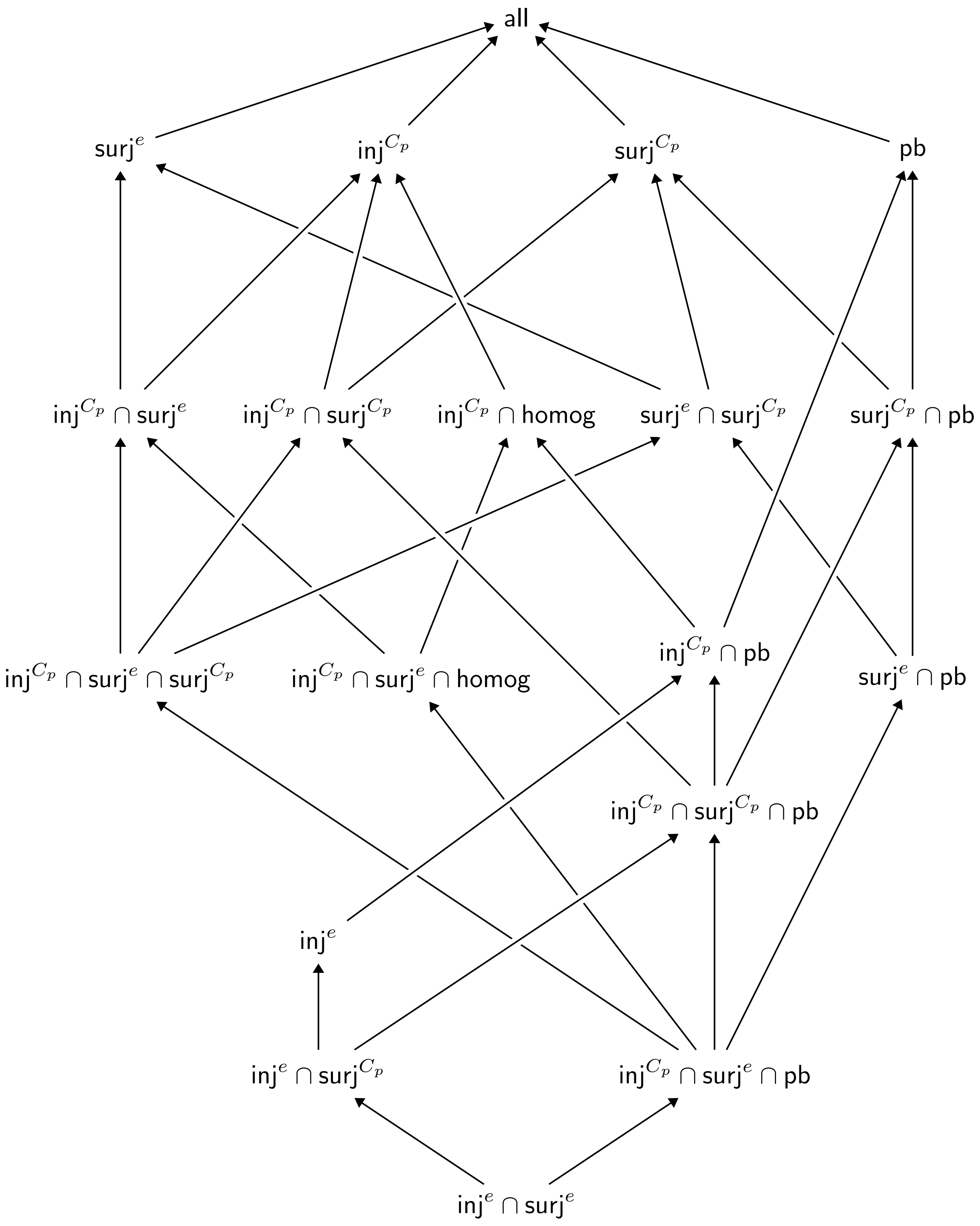
\begin{figure}[ph]
	\begin{tikzpicture}[scale=3.5]
	\node (v0) at (0, 2) {$ \all $};
	\node (v1) at (-0.5, 1.5) {$ \inj^{C_p} $};
	\node (v2) at (-1.5, 1.5) {$ \surj^e $};
	\node (v3) at (0.5, 1.5) {$ \surj^{C_p} $};
	\node (v4) at (-1.5, 0.5) {$ \inj^{C_p} \cap \surj^e $};
	\node (v5) at (-0.75, 0.5) {$ \inj^{C_p} \cap \surj^{C_p} $};
	\node (v6) at (0, 0.5) {$ \inj^{C_p} \cap \homog $};
	\node (v7) at (0.75, 0.5) {$ \surj^e \cap \surj^{C_p} $};
	\node (v8) at (1.5, 1.5) {$ \pb $};
	\node (v9) at (-1.5, -0.5) {$ \inj^{C_p} \cap \surj^e \cap \surj^{C_p} $};
	\node (v10) at (-0.4, -0.5) {$ \inj^{C_p} \cap \surj^e \cap \homog $};
	\node (v11) at (0.75, -0.4) {$ \inj^{C_p} \cap \pb $};
	\node (v12) at (1.5, 0.5) {$ \surj^{C_p} \cap \pb $};
	\node (v13) at (-0.75, -1.5) {$ \inj^e $};
	\node (v14) at (0.75, -1) {$ \inj^{C_p} \cap \surj^{C_p} \cap \pb $};
	\node (v15) at (1.5, -0.5) {$ \surj^e \cap \pb $};
	\node (v16) at (-0.75, -2) {$ \inj^e \cap \surj^{C_p} $};
	\node (v17) at (0.75, -2) {$ \inj^{C_p} \cap \surj^e \cap \pb $};
	\node (v18) at (0, -2.5) {$ \inj^e \cap \surj^e $};
	\draw[white, ultra thick, double=black, arrows={-Triangle[black,scale=0.25]}] (v18) -- (v17);
	\draw[white, ultra thick, double=black, arrows={-Triangle[black,scale=0.25]}] (v18) -- (v16);
	\draw[white, ultra thick, double=black, arrows={-Triangle[black,scale=0.25]}] (v17) -- (v15);
	\draw[white, ultra thick, double=black, arrows={-Triangle[black,scale=0.25]}] (v17) -- (v10);
	\draw[white, ultra thick, double=black, arrows={-Triangle[black,scale=0.25]}] (v17) -- (v9);
	\draw[white, ultra thick, double=black, arrows={-Triangle[black,scale=0.25]}] (v17) -- (v14);
	\draw[white, ultra thick, double=black, arrows={-Triangle[black,scale=0.25]}] (v15) -- (v7);
	\draw[white, ultra thick, double=black, arrows={-Triangle[black,scale=0.25]}] (v15) -- (v12);
	\draw[white, ultra thick, double=black, arrows={-Triangle[black,scale=0.25]}] (v10) -- (v4);
	\draw[white, ultra thick, double=black, arrows={-Triangle[black,scale=0.25]}] (v10) -- (v6);
	\draw[white, ultra thick, double=black, arrows={-Triangle[black,scale=0.25]}] (v9) -- (v7);
	\draw[white, ultra thick, double=black, arrows={-Triangle[black,scale=0.25]}] (v9) -- (v4);
	\draw[white, ultra thick, double=black, arrows={-Triangle[black,scale=0.25]}] (v9) -- (v5);
	\draw[white, ultra thick, double=black, arrows={-Triangle[black,scale=0.25]}] (v7) -- (v2);
	\draw[white, ultra thick, double=black, arrows={-Triangle[black,scale=0.25]}] (v7) -- (v3);
	\draw[white, ultra thick, double=black, arrows={-Triangle[black,scale=0.25]}] (v4) -- (v2);
	\draw[white, ultra thick, double=black, arrows={-Triangle[black,scale=0.25]}] (v4) -- (v1);
	\draw[white, ultra thick, double=black, arrows={-Triangle[black,scale=0.25]}] (v2) -- (v0);
	\draw[white, ultra thick, double=black, arrows={-Triangle[black,scale=0.25]}] (v16) -- (v14);
	\draw[white, ultra thick, double=black, arrows={-Triangle[black,scale=0.25]}] (v16) -- (v13);
	\draw[white, ultra thick, double=black, arrows={-Triangle[black,scale=0.25]}] (v14) -- (v12);
	\draw[white, ultra thick, double=black, arrows={-Triangle[black,scale=0.25]}] (v14) -- (v5);
	\draw[white, ultra thick, double=black, arrows={-Triangle[black,scale=0.25]}] (v14) -- (v11);
	\draw[white, ultra thick, double=black, arrows={-Triangle[black,scale=0.25]}] (v12) -- (v3);
	\draw[white, ultra thick, double=black, arrows={-Triangle[black,scale=0.25]}] (v12) -- (v8);
	\draw[white, ultra thick, double=black, arrows={-Triangle[black,scale=0.25]}] (v5) -- (v3);
	\draw[white, ultra thick, double=black, arrows={-Triangle[black,scale=0.25]}] (v5) -- (v1);
	\draw[white, ultra thick, double=black, arrows={-Triangle[black,scale=0.25]}] (v3) -- (v0);
	\draw[white, ultra thick, double=black, arrows={-Triangle[black,scale=0.25]}] (v13) -- (v11);
	\draw[white, ultra thick, double=black, arrows={-Triangle[black,scale=0.25]}] (v11) -- (v8);
	\draw[white, ultra thick, double=black, arrows={-Triangle[black,scale=0.25]}] (v11) -- (v6);
	\draw[white, ultra thick, double=black, arrows={-Triangle[black,scale=0.25]}] (v8) -- (v0);
	\draw[white, ultra thick, double=black, arrows={-Triangle[black,scale=0.25]}] (v6) -- (v1);
	\draw[white, ultra thick, double=black, arrows={-Triangle[black,scale=0.25]}] (v1) -- (v0);
	\node[fill=white] at (v0) {$ \all $};
	\node[fill=white] at (v1) {$ \inj^{C_p} $};
	\node[fill=white] at (v2) {$ \surj^e $};
	\node[fill=white] at (v3) {$ \surj^{C_p} $};
	\node[fill=white] at (v4) {$ \inj^{C_p} \cap \surj^e $};
	\node[fill=white] at (v5) {$ \inj^{C_p} \cap \surj^{C_p} $};
	\node[fill=white] at (v6) {$ \inj^{C_p} \cap \homog $};
	\node[fill=white] at (v7) {$ \surj^e \cap \surj^{C_p} $};
	\node[fill=white] at (v8) {$ \pb $};
	\node[fill=white] at (v9) {$ \inj^{C_p} \cap \surj^e \cap \surj^{C_p} $};
	\node[fill=white] at (v10) {$ \inj^{C_p} \cap \surj^e \cap \homog $};
	\node[fill=white] at (v11) {$ \inj^{C_p} \cap \pb $};
	\node[fill=white] at (v12) {$ \surj^{C_p} \cap \pb $};
	\node[fill=white] at (v13) {$ \inj^e $};
	\node[fill=white] at (v14) {$ \inj^{C_p} \cap \surj^{C_p} \cap \pb $};
	\node[fill=white] at (v15) {$ \surj^e \cap \pb $};
	\node[fill=white] at (v16) {$ \inj^e \cap \surj^{C_p} $};
	\node[fill=white] at (v17) {$ \inj^{C_p} \cap \surj^e \cap \pb $};
	\node[fill=white] at (v18) {$ \inj^e \cap \surj^e $};
\end{tikzpicture}
	\caption{The meet-semilattice of all intersections $\bigcap_{i \in I} \rlp{\{f_i\}}$, where $f_i$ is a morphism in $\CpSet$ with transitive codomain. This meet-semilattice is generated by our code \cite{code} as a sub-meet-semilattice of \Cref{fig:primordial_meet_semilattice}. }
	\label{figure:transitive_generated_rlcs}
\end{figure}

Nevertheless, \Cref{fig:primordial_meet_semilattice} \textit{contains} all possible intersections $\bigcap_{i \in I} \rlp{\{f_i\}}$. We may extract them by taking the sub-meet-semilattice of \Cref{fig:primordial_meet_semilattice} generated by those classes which \textit{do} appear as $\rlp{\{f_i\}}$ for some morphism $f_i$ with transitive codomain, i.e.\ those classes in \cref{equation:right lifting classes}. We generate this lattice in our code \cite{code}, and display it in \Cref{figure:transitive_generated_rlcs}.

\bigskip

At this point, we understand both pieces of \cref{equation:rlp of coproduct} independently:  \Cref{figure:transitive_generated_rlcs} describes all possible options for $\bigcap_i \rlp{\{f_i\}}$, and \Cref{fig:Cp empty classes} describes all possible options for $\bigcup_i \{f_i\}^\emptyset$. To put them together, we compute all containments $A \cup A' \subseteq B \cup B'$, where $A$ and $B$ are primordial intersections (\Cref{fig:primordial_meet_semilattice}) and $A'$ and $B'$ are collateral classes (\Cref{fig:Cp empty classes}).

\begin{definition}
	A \textit{mixed class} is a union $A \cup A'$ where $A$ is a primordial intersection (\Cref{fig:primordial_meet_semilattice}) and $A'$ is a collateral class (\Cref{fig:Cp empty classes}).
\end{definition}

For this, we must first understand the containment relations between the primordial classes $\inj^e$, $\inj^{C_p}$, $\surj^e$, $\surj^{C_p}$, $\pb$, and $\homog$, and the primordial vacuous classes $\{\id_\emptyset\}$, $\{\id_\emptyset\}^{C_p}$, $\init^e$, $\init^{C_p}$, and $\init^{C_p} \cap \pb$. These relations are displayed in \cref{table:vac-in-rl containments}, which is verified in \cref{appendix:verification of table 5}.

\begin{table}[h]
	\centering
	\begin{tblr}{colspec={c|*{6}{c}}, cells={mode=math}, colsep=5pt}
	\toprule
	\downarrow \subseteq \rightarrow & \inj^e & \inj^{C_p} & \surj^e & \surj^{C_p} & \pb    & \homog \\
	\midrule
	\{\id_\emptyset\}              & \cmark & \cmark     & \cmark  & \cmark      & \cmark & \cmark \\
	\init^e                          & \cmark & \cmark     & \xmark  & \xmark      & \cmark & \cmark \\
	\{\id_\emptyset\}^{C_p}        & \xmark & \cmark     & \xmark  & \cmark      & \cmark & \cmark \\
	\init^{C_p} \cap \pb             & \xmark & \cmark     & \xmark  & \xmark      & \cmark & \cmark \\
	\init^{C_p}                      & \xmark & \cmark     & \xmark  & \xmark      & \xmark & \cmark \\
	\bottomrule
\end{tblr}

	\caption{Containments of the vacuous right-lifting classes $\{\id_\emptyset\}$, $\{\id_\emptyset\}^{C_p}$, $\init^e$, $\init^{C_p}$, and $\init^{C_p} \cap \pb$ inside the classes $\inj^e$, $\inj^{C_p}$, $\surj^e$, $\surj^{C_p}$, $\pb$, and $\homog$. The correctness of this table is verified in \cref{appendix:verification of table 5}.} 
	\label{table:vac-in-rl containments}
\end{table}

\begin{proposition}[\computerassisted]
	\label{ordering on rlp of morphisms}
	Let $A \cup A'$ and $B \cup B'$ be mixed classes. Then $A \cup A' \subseteq B \cup B'$ if and only if $A \subseteq B$ and, for every primordial vacuous class $C$ \eqref{equation:vacuous classes} with $C \subseteq A'$, we have $C \subseteq B$ or $C \subseteq B'$.
\end{proposition}

\begin{proof}
	The backwards direction is easy. To check the forwards direction, we rely on computer verification. The strategy is as follows: we define a relation on classes of the form $A \cup A'$ by
	\begin{align*}
		A \cup A' \preceq B \cup B' \iff{} & A \subseteq B \text{ and, for every class $C$ in \eqref{equation:vacuous classes} with $C \subseteq A'$,} \\
			                                   & C \subseteq B \text{ or } C \subseteq B'.
	\end{align*}
	We verify by computer that this in fact defines a partial order on these classes,\footnote{In the code, a mixed class is stored as a pair $(A, A')$. Distinct pairs may represent the same class; for instance, $(\all, A')$ represents $\all$ for every collateral class $A'$. So the code first puts each pair into a canonical form, removing from $A'$ every class in \eqref{equation:vacuous classes} that is already contained in $A$. We then check that $\preceq$ is a partial order on these canonical pairs.} and generate a poset $(P,\preceq)$ of mixed classes from this partial order. This is the part that requires \cref{table:vac-in-rl containments}. Because we've shown the backwards direction of this proposition already, we know that this relation $\preceq$ is no stronger than subset containment $\subseteq$. Therefore, the poset $(P,\preceq)$ can only differ from the true poset $(P, \subseteq)$ of containments by failing to record some true containments. Finally, we check (by computer) that for each non-containment recorded in $(P, \preceq)$ that there is an example morphism in \cref{appendix:proper containments} that witnesses the non-containment (\cref{strategy}). Therefore, $(P, \preceq)$ is exactly the same poset as $(P, \subseteq)$. This is implemented in \cite{code}.
\end{proof}

\begin{figure}[h]
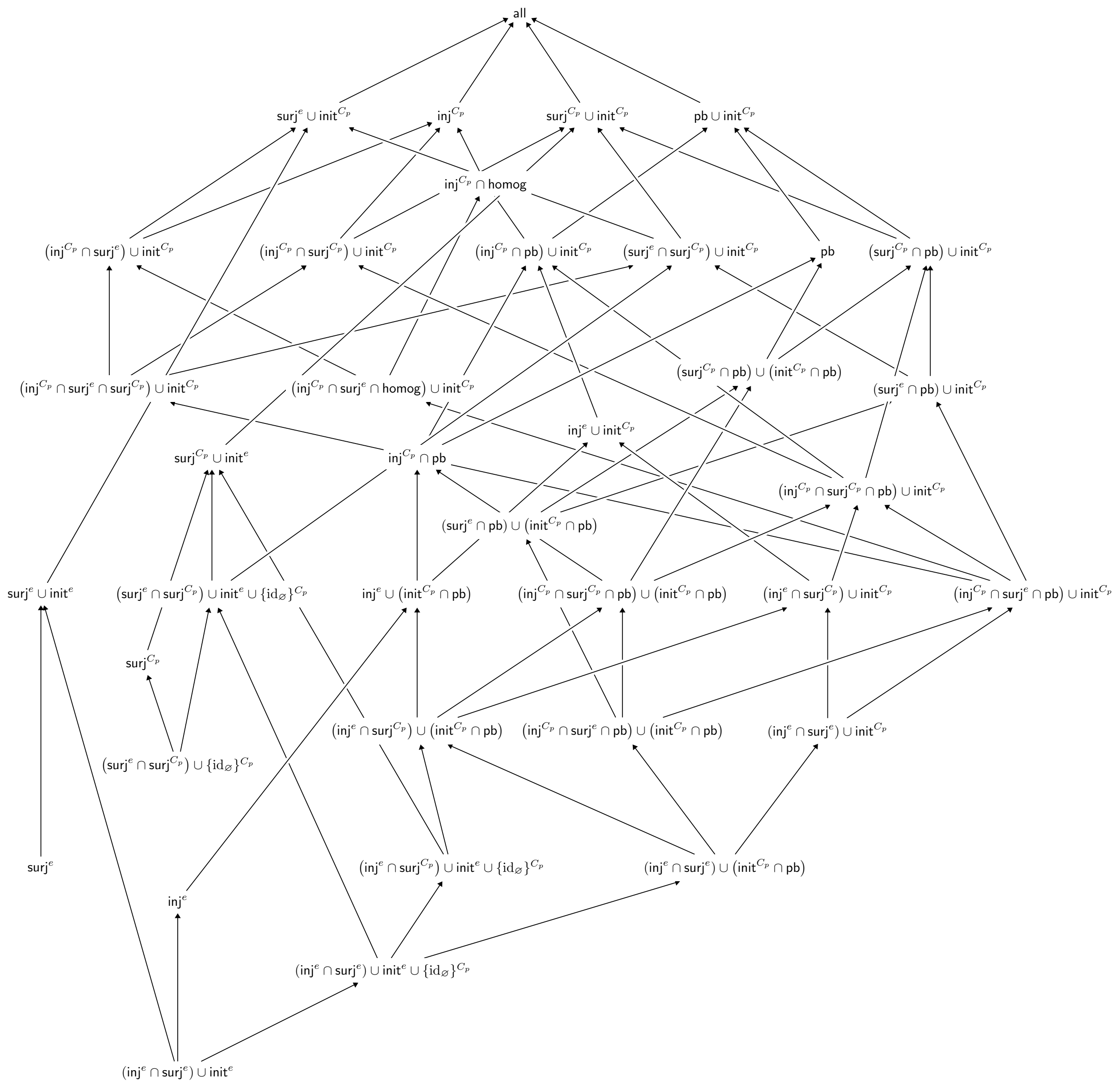

	\includestandalone[width=\textwidth]{figures/singleton_rlcs}
	\caption{The poset of right-lifting classes $\rlp{\{f\}}$, where $f$ is an arbitrary morphism of $C_p$-sets. This figure is generated by computer, and differs from the output of our code \textit{only} by changing the positions of nodes for readability. A larger version of this figure is available in the GitHub repository accompanying this paper \cite{code}.}
	\label{figure:singleton rlcs}
\end{figure}

This last proposition gives us a method to compute the poset of right-lifting classes for individual morphisms: suppose that $f = \coprod_i f_i$ and $g = \coprod_j g_j$, where $f_i$ and $g_j$ have transitive codomain. Then we have
\[
	\rlp{\{f\}} =
	\left(\bigcap_i \rlp{\{f_i\}}\right)
	\cup
	\left(\bigcup_i \{f_i\}^\emptyset\right)
	\quad \text{ and } \quad
	\rlp{\{g\}} =
	\left(\bigcap_j \rlp{\{g_j\}}\right)
	\cup
	\left(\bigcup_j \{g_j\}^\emptyset\right).
\]
Now \cref{ordering on rlp of morphisms} says that $\rlp{\{f\}} \subseteq \rlp{\{g\}}$ if and only if
\[
	\bigcap_i \rlp{\{f_i\}}
	\subseteq
	\bigcap_j \rlp{\{g_j\}}
\]
and, for every $i$,
\[
	\left(
	\{f_i\}^\emptyset
	\subseteq
	\bigcap_j \rlp{\{g_j\}}
	\quad \text{ or } \quad
	\{f_i\}^\emptyset
	\subseteq
	\bigcup_j \{g_j\}^\emptyset
	\right).
\]
Here we use that each $\{f_i\}^\emptyset$ is one of the classes in \eqref{equation:vacuous classes}, and that every class in \eqref{equation:vacuous classes} contained in $\bigcup_i \{f_i\}^\emptyset$ is contained in some $\{f_i\}^\emptyset$ (see \Cref{fig:Cp empty classes}).
We use this to implement code \cite{code} to compute the poset of all right-lifting classes, displayed in \Cref{figure:singleton rlcs}. We take this opportunity to summarize what we've done thus far to create this poset.

\begin{theorem}[\computerassisted]
	The poset of possible right-lifting classes $\rlp{\{f\}}$ of a single morphism $f$ in $\CpSet$ is shown in \Cref{figure:singleton rlcs}.
\end{theorem}

\begin{proof}
	By \cref{lifting class of coproduct}, the right-lifting class of any individual morphism is of the form
	\[
		\rlp{\{f\}} =
		\left(\bigcap_i \rlp{\{f_i\}}\right)
		\cup
		\left(\bigcup_i \{f_i\}^\emptyset\right).
	\]
	Let $A = \bigcap_i \rlp{\{f_i\}}$, and let $A' = \bigcup_i \{f_i\}^\emptyset$. All possible options for $A$ are described in \Cref{figure:transitive_generated_rlcs}, and all possible options for $A'$ are described in \Cref{fig:Cp empty classes}. We note that not every pair $(A,A')$ of classes from \Cref{figure:transitive_generated_rlcs,fig:Cp empty classes} appears as $\rlp{\{f\}} = A \cup A'$; those classes that do appear are constrained by the morphisms in $\CpSet$, insofar as we can only have $\rlp{\{f\}} = A \cup A'$ when $A = \bigcap_i \rlp{\{f_i\}}$ and $A' = \bigcup_i \{f_i\}^\emptyset$ for the same set $\{f_i\}_i$.  We take this constraint into account in our code \cite{code}.

	The containment relation for
	\[
		\rlp{\{f\}} = A \cup A' \subseteq \rlp{\{g\}} = B \cup B'
	\]
	is described in \cref{ordering on rlp of morphisms}. We implement this in code \cite{code}, and display the result in \Cref{figure:singleton rlcs}.
\end{proof}

\subsection{Arbitrary Right-Lifting Classes}

\Cref{figure:singleton rlcs} shows all possible right-lifting classes of a single morphism. To find the right-lifting classes of an arbitrary collection of morphisms, we must find the closure of the poset of singleton right-lifting classes in \Cref{figure:singleton rlcs} under intersections (\cref{lifting class of union}). \emph{A priori}, this would appear to only compute right-lifting classes of finite sets of morphisms. However, there are only finitely many distinct classes $\rlp{\{f\}}$ (\Cref{figure:singleton rlcs}), so any intersection of such classes is equal to a finite intersection.

To efficiently compute all intersections of classes in \Cref{figure:singleton rlcs}, we study containments of intersections of these classes. As in \cref{equation:rlp of coproduct}, write $\rlp{\{f\}} = A \cup A'$ where $A = \bigcap_i \rlp{\{f_i\}}$ and $A' = \bigcup_i \{f_i\}^\emptyset$, where $f = \coprod_i f_i$ and each $f_i$ has transitive codomain. Similarly, for a morphism $g = \coprod_i g_i$, write $\rlp{\{g\}} = B \cup B'$ with $B = \bigcap_i \rlp{\{g_i\}}$ and $B' = \bigcup_i \{g_i\}^\emptyset$.  Then
\[
	\rlp{\{f\}} \cap \rlp{\{g\}} = (A \cap B) \cup (A \cap B') \cup (A' \cap B) \cup (A' \cap B').
\]
Note that $A \cap B$ is a primordial intersection and already appears in \Cref{figure:transitive_generated_rlcs}, and $A' \cap B'$ is an intersection of collateral classes from \Cref{fig:Cp empty classes}. The other two terms $A \cap B'$ and $A' \cap B$ have not yet been computed. Taking the intersection of more than two right-lifting-classes of single morphisms would yield a similar decomposition
	\begin{equation}
		\label{intersection of rlc of three morphisms}
		\rlp{\{f\}} \cap \rlp{\{g\}} \cap \rlp{\{h\}} \cap \cdots = A \cup \bigcup_k Q_k,
	\end{equation}
	where $A$ is a class in \Cref{figure:transitive_generated_rlcs} and each $Q_k$ is a finite intersection of one or more collateral classes with zero or more primordial classes.

	Note that the poset of collateral classes (\Cref{fig:Cp empty classes}) is \emph{almost} closed under intersections, but there is one intersection that does not already appear in that figure. In particular, $\{\id_\emptyset\}^{C_p} \cap \init^e$ is distinct from $\{\id_\emptyset\}$, which is witnessed by the morphism $\emptyset \to C_p/e$.
	Thus, each $Q_k$ is a class that may be obtained by closing the collateral classes, plus this new intersection $\{\id_\emptyset\}^{C_p} \cap \init^e$, under intersections with the six primordial classes. We compute all possibilities for the classes $Q_k$ in \cref{table:big table} below.

\begin{definition}
	A \textit{collateral intersection} is a class $Q$ that is the intersection of one or more collateral classes with zero or more primordial classes; these are exactly the classes that appear in \cref{table:big table}.
\end{definition}

\begin{table}[h]
	\centering
	\resizebox{\textwidth}{!}{%
		\IntersectionTable{
			cell{2-4,6}{2,4} = {fg=red},
			cell{9-14}{1}    = {fg=red},
			cell{9}{4}       = {fg=red},
		}%
	}
	\vspace{1em}
	\caption{Iterated intersections of the six primordial classes $\inj^e$, $\inj^{C_p}$, $\surj^e$, $\surj^{C_p}$, $\pb$, and $\homog$ with the seven intersections of collateral classes $\init^{e}$, $\init^{C_p}$, $\init^{C_p} \cap \pb$, $\{\id_\emptyset\}^{C_p}$, $\{\id_\emptyset\}$, $\init^e \cup \{\id_\emptyset\}^{C_p}$ and $\{\id_\emptyset\}^{C_p} \cap \init^e$. We produce this table iteratively, beginning with the first seven rows. When an entry in the table is distinct from the existing row or column labels, we highlight it in red and add this entry as a new row, repeating until the computation stabilizes.
			Note that every intersection is at most a threefold intersection, and threefold intersections do implicitly appear because $\bij^e = \inj^e \cap \surj^e$. We verify the correctness of this table in \cref{appendix:table 6}. Classes in this table are called \textit{collateral intersections.}}
	\label{table:big table}
\end{table}

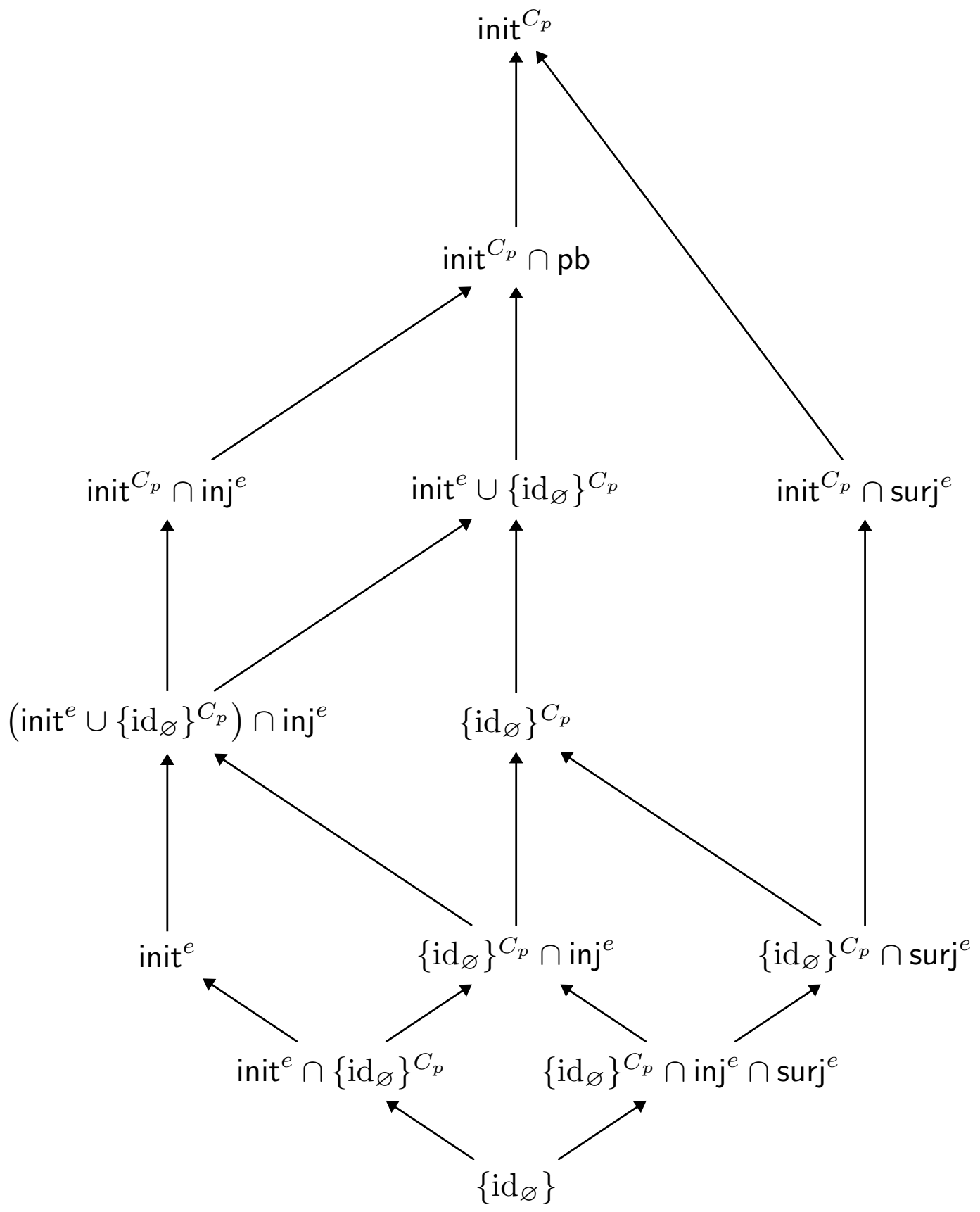
\begin{figure}[h]
	\begin{tikzpicture}[scale=2.5]
    \node (v0) at (0, -2) {$ \{\id_\emptyset\} $};
    \node (v1) at (-1.5, -1) {$ \init^e $};
    \node (v2) at (-0.75, -1.5) {$ \init^e \cap \{\id_\emptyset\}^{C_p} $};
    \node (v3) at (0, 0) {$ \{\id_\emptyset\}^{C_p} $};
    \node (v4) at (0, 1) {$ \init^e \cup \{\id_\emptyset\}^{C_p} $};
    \node (v5) at (0, 2) {$ \init^{C_p} \cap \pb $};
    \node (v6) at (0, 3) {$ \init^{C_p} $};
    \node (v7) at (1.5, -1) {$ \{\id_\emptyset\}^{C_p} \cap \surj^e $};
    \node (v8) at (1.5, 1) {$ \init^{C_p} \cap \surj^e $};
    \node (v9) at (0, -1) {$ \{\id_\emptyset\}^{C_p} \cap \inj^e $};
    \node (v10) at (-1.5, 0) {$ \left(\init^e \cup \{\id_\emptyset\}^{C_p}\right) \cap \inj^e $};
    \node (v11) at (-1.5, 1) {$ \init^{C_p} \cap \inj^e $};
    \node (v12) at (0.75, -1.5) {$ \{\id_\emptyset\}^{C_p} \cap \inj^e \cap \surj^e $};
    \draw[white, ultra thick, double=black, arrows={-Triangle[black,scale=0.25]}] (v0) -- (v2);
    \draw[white, ultra thick, double=black, arrows={-Triangle[black,scale=0.25]}] (v0) -- (v12);
    \draw[white, ultra thick, double=black, arrows={-Triangle[black,scale=0.25]}] (v2) -- (v1);
    \draw[white, ultra thick, double=black, arrows={-Triangle[black,scale=0.25]}] (v2) -- (v9);
    \draw[white, ultra thick, double=black, arrows={-Triangle[black,scale=0.25]}] (v1) -- (v10);
    \draw[white, ultra thick, double=black, arrows={-Triangle[black,scale=0.25]}] (v12) -- (v7);
    \draw[white, ultra thick, double=black, arrows={-Triangle[black,scale=0.25]}] (v12) -- (v9);
    \draw[white, ultra thick, double=black, arrows={-Triangle[black,scale=0.25]}] (v7) -- (v8);
    \draw[white, ultra thick, double=black, arrows={-Triangle[black,scale=0.25]}] (v7) -- (v3);
    \draw[white, ultra thick, double=black, arrows={-Triangle[black,scale=0.25]}] (v8) -- (v6);
    \draw[white, ultra thick, double=black, arrows={-Triangle[black,scale=0.25]}] (v9) -- (v3);
    \draw[white, ultra thick, double=black, arrows={-Triangle[black,scale=0.25]}] (v9) -- (v10);
    \draw[white, ultra thick, double=black, arrows={-Triangle[black,scale=0.25]}] (v3) -- (v4);
    \draw[white, ultra thick, double=black, arrows={-Triangle[black,scale=0.25]}] (v10) -- (v4);
    \draw[white, ultra thick, double=black, arrows={-Triangle[black,scale=0.25]}] (v10) -- (v11);
    \draw[white, ultra thick, double=black, arrows={-Triangle[black,scale=0.25]}] (v4) -- (v5);
    \draw[white, ultra thick, double=black, arrows={-Triangle[black,scale=0.25]}] (v11) -- (v5);
    \draw[white, ultra thick, double=black, arrows={-Triangle[black,scale=0.25]}] (v5) -- (v6);
    \node[fill=white] at (v0) {$ \{\id_\emptyset\} $};
    \node[fill=white] at (v1) {$ \init^e $};
    \node[fill=white] at (v2) {$ \init^e \cap \{\id_\emptyset\}^{C_p} $};
    \node[fill=white] at (v3) {$ \{\id_\emptyset\}^{C_p} $};
    \node[fill=white] at (v4) {$ \init^e \cup \{\id_\emptyset\}^{C_p} $};
    \node[fill=white] at (v5) {$ \init^{C_p} \cap \pb $};
    \node[fill=white] at (v6) {$ \init^{C_p} $};
    \node[fill=white] at (v7) {$ \{\id_\emptyset\}^{C_p} \cap \surj^e $};
    \node[fill=white] at (v8) {$ \init^{C_p} \cap \surj^e $};
    \node[fill=white] at (v9) {$ \{\id_\emptyset\}^{C_p} \cap \inj^e $};
    \node[fill=white] at (v10) {$ \left(\init^e \cup \{\id_\emptyset\}^{C_p}\right) \cap \inj^e $};
    \node[fill=white] at (v11) {$ \init^{C_p} \cap \inj^e $};
    \node[fill=white] at (v12) {$ \{\id_\emptyset\}^{C_p} \cap \inj^e \cap \surj^e $};
\end{tikzpicture}
	\caption{The meet-semilattice of collateral intersections, created from  \cref{table:big table}. Note that this poset contains \Cref{fig:Cp empty classes} as a sub-poset.}
	\label{fig:small lattice}
\end{figure}

\vspace*{-3mm}

From \cref{table:big table}, we can produce a meet-semilattice, with meet given by intersection, of collateral intersections. This is the meet-semilattice of all classes that can appear as $Q_k$ in \cref{intersection of rlc of three morphisms} above. We display this meet-semilattice in \Cref{fig:small lattice}. The poset structure is easily determined from the intersection structure: $\mathcal P \subseteq \mathcal Q$ if and only if $\mathcal P \cap \mathcal Q = \mathcal P$. This meet-semilattice is in fact a lattice, because it is finite and has a maximum element.

\begin{proposition}[\computerassisted]
	\label{ordering on mixed intersections}
	Let $A$ and $B$ be primordial intersections (\Cref{fig:primordial_meet_semilattice}) and let $\{P_i\}_i$ and $\{Q_j\}_j$ be collateral intersections. Then
	\[
		A \cup \bigcup_i P_i \subseteq B \cup \bigcup_j Q_j
	\]
	if and only if $A \subseteq B$ and for all $i$, $P_i \subseteq B \cup \bigcup_j Q_j$.
\end{proposition}

\begin{proof}
	The backwards direction is easy. To prove the forwards direction, we rely on computer verification. In order to use \cref{strategy} here, we must be able to have the computer tell us whether or not $P_i \subseteq B \cup \bigcup_j Q_j$. For this, we use the same intersection trick as above:
	\begin{align*}
		P_i \subseteq B \cup \bigcup_j Q_j\iff & P_i \cap (B \cup \bigcup_j Q_j) = P_i            \\
		\iff                                   & (P_i \cap B) \cup (P_i \cap \bigcup_j Q_j) = P_i \\
		\iff                                   & (P_i \cap B) \cup \bigcup_j (P_i\cap Q_j) = P_i
	\end{align*}
	Note that $P_i \cap B$ and $P_i \cap Q_j$ are all collateral intersection classes appearing in \cref{table:big table}.

	To evaluate such statements, we must understand when finite unions of collateral intersection classes are equal to a single collateral intersection class. That is, for a collection $\{R_k\}_k$ of classes $R_k$ from \Cref{table:big table}, we ask when
	\[
		\bigcup_k R_k = P
	\]
	for some $P$ a single collateral intersection class.

	We note that there is a unique candidate for $P$, because the poset \Cref{fig:small lattice} is a finite lattice with containment given by $\subseteq$. This candidate for $P$ is the join $\bigvee_k R_k$ of the classes $R_k$ in the lattice of \Cref{fig:small lattice}, which is easily read off from the figure. We also note that we always have \[\bigcup_k R_k \subseteq \bigvee_k R_k,\] but it may be the case that this containment is proper.

	So we employ \cref{strategy} in our code \cite{code}, which shows that this containment is in fact proper most of the time, with exceptions generated by three cases. The first case is that the union of $\init^e$ and $\{\id_\emptyset\}^{C_p}$ is the class $\init^e \cup \{\id_\emptyset\}^{C_p}$, which is also a collateral intersection class. The other two are
	\begin{align*}
		\left(\init^e \cup \{\id_\emptyset\}^{C_p}\right) \cap \inj^e & = \init^e \cup \left(\{\id_\emptyset\}^{C_p} \cap \inj^e\right),                                           \\
		\init^e \cup \{\id_\emptyset\}^{C_p}                          & = \{\id_\emptyset\}^{C_p} \cup \left(\left(\init^e \cup \{\id_\emptyset\}^{C_p}\right) \cap \inj^e\right).
	\end{align*}
	Together with unions of these equations with other collateral intersection classes, this describes all possible cases when a union of collateral intersection classes $\bigcup_k R_k$ is equal to a single collateral intersection class.

	Understanding this, we can now implement the relation $\preceq$ by
	\[
		A \cup \bigcup_i P_i  \preceq B \cup \bigcup_j Q_j \iff A \subseteq B \text{ and for all $i$ } P_i \subseteq B \cup \bigcup_j Q_j.
	\]
	Again, we use the computer to check that this is a well-defined partial order on all classes of this form.\footnote{As with mixed classes, the code stores such a class as a pair consisting of $A$ and a set of collateral intersections, and distinct pairs may represent the same class. So the code first puts each pair into a canonical form: it replaces the set of collateral intersections by the maximal collection of collateral intersections that are contained in $A \cup \bigcup_i P_i$ but not in $A$. It then checks that $\preceq$ is a partial order on these canonical pairs.} By the backward direction of this proposition, we also know that this relation is no stronger than subset containment $\subseteq$. Finally, we use \cref{strategy} to show that this relation $\preceq$ is in fact equal to subset containment $\subseteq$.
\end{proof}

\begin{sidewaysfigure}[p]
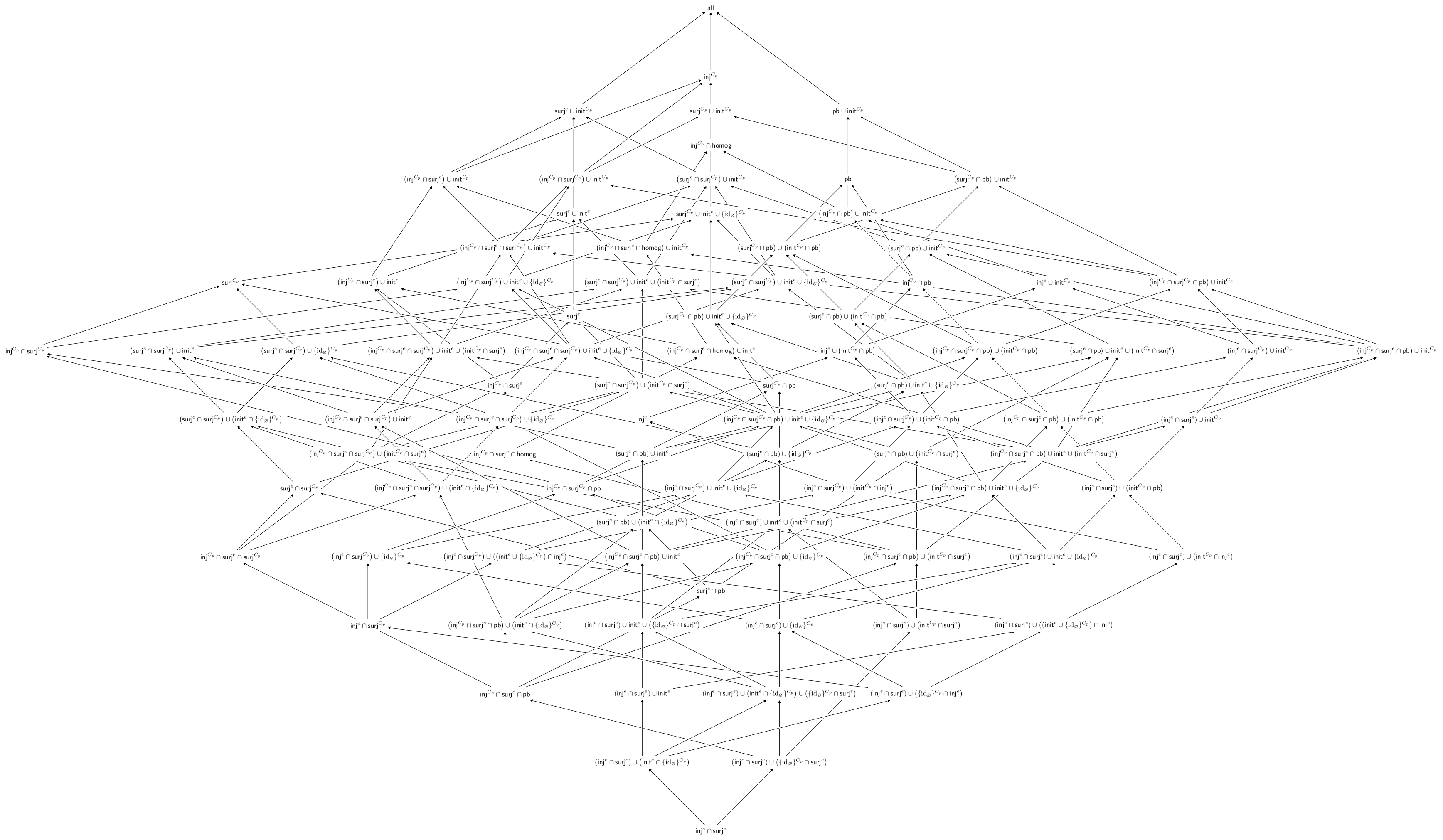

	\centering
	\vspace*{10cm}\includestandalone[width=\textwidth]{figures/all_rlc_semilattice}
	\caption{The poset of \emph{all} right-lifting classes $\mathcal{M}^\boxslash$, where $\mathcal{M}$ is an arbitrary class of morphisms in $\CpSet$. A larger version of this figure is available in the GitHub repository accompanying this paper \cite{code}.}
	\label{fig:biggest lattice}
\end{sidewaysfigure}

\Cref{ordering on mixed intersections} allows us to produce by code \cite{code} the meet-semilattice of all right lifting classes in the category $\CpSet$, shown in \Cref{fig:biggest lattice}. We again summarize what we've done in a theorem.

\begin{theorem}[\computerassisted]
	\label{thm:wfs lattice}
	The lattice of weak factorization systems on the category of $C_p$-sets is isomorphic to the lattice displayed in \Cref{fig:biggest lattice}. There are 89 weak factorization systems on $C_p$-sets.
\end{theorem}

\begin{proof}
	The lattice displayed is the lattice of all possible right-lifting classes of classes of morphisms in $C_p$-sets. Any weak factorization system $(\cat L, \cat R)$ has $\cat R = \rlp{\cat L}$, so $\cat R$ is one of the classes in \Cref{fig:biggest lattice}. Conversely, let $\cat R = \rlp{\mathcal M}$ be one of the classes in \Cref{fig:biggest lattice}, where $\mathcal M$ is a class of morphisms. By \cref{lifting class of union}, $\cat R = \bigcap_{f \in \mathcal M} \rlp{\{f\}}$. There are only finitely many distinct classes $\rlp{\{f\}}$ (\Cref{figure:singleton rlcs}), so choosing one morphism $f \in \mathcal M$ for each distinct class $\rlp{\{f\}}$ that occurs gives a finite set $I \subseteq \mathcal M$ with $\cat R = \rlp{I}$. Now by \cref{smallobjectargument}, $(\llp{\cat R}, \cat R)$ is a weak factorization system. These two assignments are inverse to each other, and they are order-preserving when weak factorization systems are ordered by $(\cat L, \cat R) \leq (\cat L', \cat R')$ if and only if $\cat R \subseteq \cat R'$. Thus, the lattice of right-lifting classes is isomorphic to that of weak factorization systems, which is \Cref{fig:biggest lattice}.
\end{proof}

\begin{corollary}[\computerassisted]
	\label{cor:final count}
	There are 1777 pre-model structures on the category of $C_p$-sets.
\end{corollary}

\begin{proof}
	Pre-model structures are equivalent to intervals (pairs $(x,y)$ with $x \leq y$) in the lattice of weak factorization systems. We count such intervals by computer \cite{code}.
\end{proof}

We note that this number $1777$ is an upper bound on the number of model structures on $\CpSet$. However, we suspect this bound is very far from being tight. Evidence for this suspicion is provided by the cases of model structures on $\Set$, where there are 18 pre-model structures by only 9 model structures, and $\textsf{Vect}_k$, where there are 9 pre-model structures but only 4 model structures \cite{ACB2025}.

Determining which of these pre-model structures are model structures is a significant undertaking that we do not attempt. For one, we have only identified a poset isomorphic to the poset of pre-model structures, which does not even completely describe each pre-model structure. To do so, we would need to compute the left-lifting classes of each class of morphisms in \Cref{fig:biggest lattice}, which is a daunting undertaking even with the aid of a computer. Then, given a pre-model structure $(\mathcal L, \mathcal R)$ and $(\mathcal L', \mathcal R')$ with $\mathcal R \subseteq \mathcal R'$, we would need to check that the class $\mathcal W = \mathcal R \circ \mathcal L'$ satisfies the two-out-of-three property. We see no way to automate these checks. Given that there are $1777$ pre-model structures on $\CpSet$, this promises to be tedious as well as difficult. Therefore, we leave the question of model structures on $\CpSet$ open. 

\clearpage
\appendix
\crefalias{section}{appendix}

\section{Table of Example Morphisms Witnessing Non-Containments}
\label{appendix:proper containments}

Recall from \cref{def quotient and fold} the definition of fold morphism
\[
	\nabla_{H,n} \colon \coprod_{n} G/H \to G/H
\]
and the quotient morphism
\[
	q^H_K \colon G/K \to G/H.
\]
Recall also the classes of morphisms in $\GSet$ from \cref{table:morphisms}.
We collect in \cref{table:primordial proper containments,table:primordial vacuous proper containments} some morphisms in $\CpSet$ built from coproducts and composition of quotients and folds, and describe which of the classes from \cref{table:morphisms} these morphisms belong to. As in \cref{section: Cp-sets}, we abbreviate $\pb^{C_p}_e$ to $\pb$. The interested reader reconstructing these tables for themselves may find \cref{easier pb and homog for Cp} helpful.

The examples in these tables are used to verify the correctness of various figures throughout \cref{section: Cp-sets}, as described in \Cref{strategy}. In order to verify a claimed non-containment $\mathcal{P} \nsubseteq \mathcal{Q}$, we must identify a morphism $f$ such that $f \in \mathcal{P}$ and $f \notin \mathcal{Q}$. For example, the non-containment $\inj^{C_p} \cap \surj^e \cap \pb \nsubseteq \inj^e$ is witnessed by the morphism $\nabla_{e,2} \colon C_p/e \amalg C_p/e \to C_p/e$ (Example 2 below). Example 14 also witnesses the same non-containment.

\begin{table}[h!]
	\centering
	\begin{tblr}{colspec={cl|*{6}{c}}, cells={mode=math}, colsep=2pt, rowsep=2pt}
		\toprule
		\text{Example} & \text{Morphism}                                        & \inj^e & \inj^{C_p} & \surj^e & \surj^{C_p} & \pb    & \homog \\
		\midrule
		0              & \nabla_{e,0}                                           & \cmark & \cmark     & \xmark  & \cmark      & \cmark & \cmark \\
		1              & \nabla_{e,1}                                           & \cmark & \cmark     & \cmark  & \cmark      & \cmark & \cmark \\
		2              & \nabla_{e,2}                                           & \xmark & \cmark     & \cmark  & \cmark      & \cmark & \cmark \\
		3              & \nabla_{C_p,0}                                         & \cmark & \cmark     & \xmark  & \xmark      & \cmark & \cmark \\
		4              & \nabla_{C_p,2}                                         & \xmark & \xmark     & \cmark  & \cmark      & \cmark & \cmark \\
		5              & q^{C_p}_e                                              & \xmark & \cmark     & \cmark  & \xmark      & \xmark & \cmark \\
		6              & \nabla_{e,1} \amalg \nabla_{e,0}                       & \cmark & \cmark     & \xmark  & \cmark      & \cmark & \cmark \\
		7              & \nabla_{C_p,2} \circ (q^{C_p}_e \amalg \nabla_{C_p,1}) & \xmark & \cmark     & \cmark  & \cmark      & \xmark & \xmark \\
		8              & \nabla_{e,1} \amalg \nabla_{C_p,0}                     & \cmark & \cmark     & \xmark  & \xmark      & \cmark & \cmark \\
		9              & \nabla_{e,2} \amalg \nabla_{C_p,0}                     & \xmark & \cmark     & \xmark  & \xmark      & \cmark & \cmark \\
		10             & q^{C_p}_e \amalg \nabla_{C_p,0}                        & \xmark & \cmark     & \xmark  & \xmark      & \xmark & \cmark \\
		11             & \nabla_{C_p,1} \amalg \nabla_{C_p,0}                   & \cmark & \cmark     & \xmark  & \xmark      & \cmark & \cmark \\
		12             & q^{C_p}_e \amalg \nabla_{C_p,1}                        & \xmark & \cmark     & \cmark  & \xmark      & \xmark & \cmark \\
		13             & \nabla_{e,1} \amalg \nabla_{e,0} \amalg \nabla_{C_p,1} & \cmark & \cmark     & \xmark  & \cmark      & \cmark & \cmark \\
		14             & \nabla_{e,2} \amalg \nabla_{C_p,1}                     & \xmark & \cmark     & \cmark  & \cmark      & \cmark & \cmark \\
		15             & \nabla_{e,2} \amalg \nabla_{e,0}                       & \xmark & \cmark     & \xmark  & \cmark      & \cmark & \cmark \\
		\bottomrule
	\end{tblr}
	\caption{Which of the primordial classes $\inj^e$, $\inj^{C_p}$, $\surj^e$, $\surj^{C_p}$, $\pb$, and $\homog$ contain each example morphism.}
	\label{table:primordial proper containments}
\end{table}

\begin{table}[h!]
	\centering
	\begin{tblr}{colspec={cl|*{5}{c}}, cells={mode=math}, colsep=2pt, rowsep=2pt}
		\toprule
		\text{Example} & \text{Morphism}                                        & \{\id_\emptyset\} & \init^e & \{\id_\emptyset\}^{C_p} & \init^{C_p} \cap \pb & \init^{C_p} \\
		\midrule
		0              & \nabla_{e,0}                                           & \xmark            & \cmark  & \cmark                  & \cmark               & \cmark      \\
		1              & \nabla_{e,1}                                           & \xmark            & \xmark  & \cmark                  & \cmark               & \cmark      \\
		2              & \nabla_{e,2}                                           & \xmark            & \xmark  & \cmark                  & \cmark               & \cmark      \\
		3              & \nabla_{C_p,0}                                         & \xmark            & \cmark  & \xmark                  & \cmark               & \cmark      \\
		4              & \nabla_{C_p,2}                                         & \xmark            & \xmark  & \xmark                  & \xmark               & \xmark      \\
		5              & q^{C_p}_e                                              & \xmark            & \xmark  & \xmark                  & \xmark               & \cmark      \\
		6              & \nabla_{e,1} \amalg \nabla_{e,0}                       & \xmark            & \xmark  & \cmark                  & \cmark               & \cmark      \\
		7              & \nabla_{C_p,2} \circ (q^{C_p}_e \amalg \nabla_{C_p,1}) & \xmark            & \xmark  & \xmark                  & \xmark               & \xmark      \\
		8              & \nabla_{e,1} \amalg \nabla_{C_p,0}                     & \xmark            & \xmark  & \xmark                  & \cmark               & \cmark      \\
		9              & \nabla_{e,2} \amalg \nabla_{C_p,0}                     & \xmark            & \xmark  & \xmark                  & \cmark               & \cmark      \\
		10             & q^{C_p}_e \amalg \nabla_{C_p,0}                        & \xmark            & \xmark  & \xmark                  & \xmark               & \cmark      \\
		11             & \nabla_{C_p,1} \amalg \nabla_{C_p,0}                   & \xmark            & \xmark  & \xmark                  & \xmark               & \xmark      \\
		12             & q^{C_p}_e \amalg \nabla_{C_p,1}                        & \xmark            & \xmark  & \xmark                  & \xmark               & \xmark      \\
		13             & \nabla_{e,1} \amalg \nabla_{e,0} \amalg \nabla_{C_p,1} & \xmark            & \xmark  & \xmark                  & \xmark               & \xmark      \\
		14             & \nabla_{e,2} \amalg \nabla_{C_p,1}                     & \xmark            & \xmark  & \xmark                  & \xmark               & \xmark      \\
		15             & \nabla_{e,2} \amalg \nabla_{e,0}                       & \xmark            & \xmark  & \cmark                  & \cmark               & \cmark      \\
		\bottomrule
	\end{tblr}
	\caption{Which of the vacuous classes $\{\id_\emptyset\}$, $\init^e$, $\{\id_\emptyset\}^{C_p}$, $\init^{C_p} \cap \pb$, and $\init^{C_p}$ contain each example morphism.}
	\label{table:primordial vacuous proper containments}
\end{table}

\newpage

\section{Verification of \texorpdfstring{\Cref{table:vac-in-rl containments}}{Table \getrefnumber{table:vac-in-rl containments}}}
\label{appendix:verification of table 5}

Here we verify the correctness of \Cref{table:vac-in-rl containments}. We replicate the table here for convenience.
\begin{center}
	\begin{tblr}{colspec={c|*{6}{c}}, cells={mode=math}, colsep=5pt}
	\toprule
	\downarrow \subseteq \rightarrow & \inj^e & \inj^{C_p} & \surj^e & \surj^{C_p} & \pb    & \homog \\
	\midrule
	\{\id_\emptyset\}              & \cmark & \cmark     & \cmark  & \cmark      & \cmark & \cmark \\
	\init^e                          & \cmark & \cmark     & \xmark  & \xmark      & \cmark & \cmark \\
	\{\id_\emptyset\}^{C_p}        & \xmark & \cmark     & \xmark  & \cmark      & \cmark & \cmark \\
	\init^{C_p} \cap \pb             & \xmark & \cmark     & \xmark  & \xmark      & \cmark & \cmark \\
	\init^{C_p}                      & \xmark & \cmark     & \xmark  & \xmark      & \xmark & \cmark \\
	\bottomrule
\end{tblr}

\end{center}
Every indicated non-containment ($\xmark$) in the table is witnessed by an example morphism from \Cref{appendix:proper containments}, which we check in our code~\cite{code}. We note that some containments ($\cmark$) imply others: $A \subseteq B$ implies $A' \subseteq B'$ for any $A' \subseteq A$ and $B' \supseteq B$. We also make use of \cref{easier pb and homog for Cp}. We now use these ideas to verify all of the claimed containments:

\begin{description}
	\item[Row 5] The containment $\init^{C_p} \subseteq \inj^{C_p}$ follows from the fact that a morphism with empty domain is injective. To verify the $\init^{C_p} \subseteq \homog$, we use \cref{easier pb and homog for Cp}. Let $f \colon X \to Y$ be a morphism in $\init^{C_p}$. Note that because $f \in \init^{C_p}$, $X$ cannot have any $C_p$-fixed points. Then for any $y \in Y$, $f^{-1}(\{y\})$ must be disjoint from $X^{C_p}$, so $f \in \homog$ by the lemma.
	      \smallskip
	\item[Row 4] The containments $\init^{C_p} \cap \pb \subseteq \inj^{C_p}$ and $\init^{C_p} \cap \pb \subseteq \homog$ follow from the verification of Row 5, because $\init^{C_p} \cap \pb \subseteq \init^{C_p}$.
	      The containment $\init^{C_p} \cap \pb \subseteq \pb$ is clear.
	      \smallskip
	\item[Row 3] The containments $\{\id_\emptyset\}^{C_p} \subseteq \inj^{C_p}$ and $\{\id_\emptyset\}^{C_p} \subseteq \surj^{C_p}$ follow from the fact $\id_\emptyset$ is both injective and surjective. The containment $\{\id_\emptyset\}^{C_p} \subseteq \pb$ follows from \cref{easier pb and homog for Cp}. Since $\pb \subseteq \homog$, then this completes the verification of Row 3.
	      \smallskip
	\item[Row 2] The containment $\init^e \subseteq \inj^e$ says simply that morphisms from $\emptyset$ are always injective. The containment $\init^e \subseteq \inj^{C_p}$ then follows from the fact $\inj^e \subseteq \inj^{C_p}$. Likewise, $\inj^e \subseteq \pb$ and $\inj^e \subseteq \homog$, so this completes the verification of the second row of the table.
	      \smallskip
	\item[Row 1] The first row of the table consists entirely of $\cmark$s because $\id_\emptyset$ is an isomorphism, and all primordial classes contain the class of isomorphisms.
\end{description}

\section{Verification of \texorpdfstring{\Cref{table:big table}}{Table \getrefnumber{table:big table}}}
\label{appendix:table 6}

Here we verify the correctness of \Cref{table:big table}. We include the table here again for convenience, with colors to be explained shortly:
\begin{center}
	\resizebox{\textwidth}{!}{%
		\IntersectionTable{%
			cell{2-14}{2-7}     = {fg=firered},
			cell{2}{3,7}        = {fg=frenchblue},
			cell{3}{3,6,7}      = {fg=frenchblue},
			cell{4}{3,6,7}      = {fg=frenchblue},
			cell{5}{2,3,6,7}    = {fg=frenchblue},
			cell{6}{3,5-7}      = {fg=frenchblue},
			cell{7}{2,3,5-7}    = {fg=frenchblue},
			cell{8}{3,2-7}      = {fg=frenchblue},
			cell{9}{2,3,6,7}    = {fg=frenchblue},
			cell{10}{2,3,6,7}   = {fg=frenchblue},
			cell{11}{3,4,7}     = {fg=frenchblue},
			cell{12}{2,3,5-7}   = {fg=frenchblue},
			cell{13}{3-7}	    = {fg=frenchblue},
			cell{14}{2-7}	    = {fg=frenchblue},
			cell{2}{2,4,6}      = {fg=forestgreen},
			cell{4}{2}          = {fg=forestgreen},
			cell{6}{2,4}	    = {fg=forestgreen},
			cell{12}{4}		    = {fg=forestgreen},
			cell{13}{2}		    = {fg=forestgreen},
			cell{2}{5}          = {appto={~(\ref{initCp cap surjCp})}},
			cell{3}{2}		    = {appto={~(\ref{(initCp cap pb) cap inje})}},
			cell{3}{4}		    = {appto={~(\ref{(initCp cap pb) cap surje})}},
			cell{3}{5}		    = {appto={~(\ref{(initCp cap pb) cap surjCp})}},
			cell{4}{4}		    = {appto={~(\ref{(inite cup idemptyCp) cap surje})}},
			cell{4}{5}		    = {appto={~(\ref{(inite cup idemptyCp) cap surjCp})}},
			cell{5}{4}		    = {appto={~(\ref{inite cap surje})}},
			cell{5}{5}		    = {appto={~(\ref{inite cap surjCp})}},
			cell{7}{4}		    = {appto={~(\ref{(idemptyCp cap inite) cap surje})}},
			cell{9}{4}		    = {appto={~(\ref{(initCp cap inje) cap surje})}},
			cell{9}{5}		    = {appto={~(\ref{(initCp cap inje) cap surjCp})}},
			cell{10}{4}		    = {appto={~(\ref{((inite cup idemptyCp) cap inje) cap surje})}},
			cell{10}{5}		    = {appto={~(\ref{((inite cup idemptyCp) cap inje) cap surjCp})}},
			cell{11}{2}		    = {appto={~(\ref{(initCp cap surje) cap inje})}},
			cell{11}{5}		    = {appto={~(\ref{(initCp cap surje) cap surjCp})}},
			cell{11}{6}		    = {appto={~(\ref{(initCp cap surje) cap pb})}},
		}
	}
\end{center}
We described how this table was produced in the caption of \Cref{table:big table}. To verify that this table is correct, we must perform the following steps:
\begin{enumerate}[(1)]
	\item Programmatically verify that all thirteen row labels represent distinct classes, using \cref{strategy} in our code \cite{code}.
	\item Each cell of the table represents an intersection $P \cap Q$, where $P$ is the class labeling the row and $Q$ is the class labeling the column. There are three cases to consider:
	      \begin{enumerate}[(a)]
		      \item The entry in the cell is either $P$ or $Q$, in which case the cell witnesses a containment $P \subseteq Q$ or $Q \subseteq P$, which we can verify programmatically using prior results. This relies, for example, on the understanding of the poset of primordial classes in \Cref{fig:primordial_meet_semilattice} and on the containments between vacuous classes and primordial classes in \cref{table:vac-in-rl containments}. These cells are highlighted blue above.
		      \item The entry in the cell is literally written as $P \cap Q$, in which case the cell is certainly correct. These cells are highlighted green above. Recall that we write $\bij^e$ as shorthand for $\inj^e \cap \surj^e$.
		      \item The entry in the cell is a new class $R$ distinct from both $P$ and $Q$. These cells are highlighted red above. In this case, we provide a proof that $P \cap Q = R$, as one of the lemmas below. The forward references in the table point to these lemmas.
	      \end{enumerate}
\end{enumerate}

\begin{lemma}
	\label{(initCp cap pb) cap inje}
	$\left(\init^{C_p} \cap \pb\right) \cap \inj^e = \init^{C_p} \cap \inj^e$
\end{lemma}
\begin{proof}
	Follows from the fact that $\inj^e \subseteq \pb$.
\end{proof}

\begin{lemma}
	\label{inite cap surje}
	$\init^e \cap \surj^e = \{\id_\emptyset\}$
\end{lemma}

\begin{proof}
	A surjective function with domain $\emptyset$ can only be $\id_\emptyset$.
\end{proof}
\begin{lemma}
	\label{inite cap surjCp}
	$\init^e \cap \surj^{C_p} = \{\id_\emptyset\}^{C_p} \cap \init^e$
\end{lemma}

\begin{proof}
	Let $f \in \init^e \cap \surj^{C_p}$. Then the domain of $f$ is $\emptyset$, so certainly the domain of $f^{C_p}$ is also $\emptyset$. Since $f \in \surj^{C_p}$, $f^{C_p}$ is a surjective function with domain $\emptyset$, so $f^{C_p} = \id_\emptyset$. Hence, $f \in \{\id_\emptyset\}^{C_p}$. Since $f \in \init^e$ by assumption, then $f \in \{\id_\emptyset\}^{C_p} \cap \init^e$.

	Conversely, we have $\{\id_\emptyset\}^{C_p} \subseteq \surj^{C_p}$, so taking intersections with $\init^e$ gives the desired inclusion.
\end{proof}
\begin{lemma}
	\label{initCp cap surjCp}
	$\init^{C_p} \cap \surj^{C_p} = \{\id_\emptyset\}^{C_p}$
\end{lemma}

\begin{proof}
	If $f \in \init^{C_p} \cap \surj^{C_p}$, then $f^{C_p}$ is a surjective function with domain $\emptyset$. In other words, $f^{C_p} = \id_\emptyset$.
\end{proof}
\begin{lemma}
	\label{(idemptyCp cap inite) cap surje}
	$\left(\{\id_\emptyset\}^{C_p} \cap \init^e\right) \cap \surj^e = \{\id_\emptyset\}$
\end{lemma}

\begin{proof}
	This follows from \cref{inite cap surje}.
\end{proof}
\begin{lemma}
	\label{(initCp cap pb) cap surjCp}
	$\left(\init^{C_p} \cap \pb \right) \cap \surj^{C_p} = \{\id_\emptyset\}^{C_p}$
\end{lemma}

\begin{proof}
	By \cref{initCp cap surjCp}, we have
	\(
	\left(\init^{C_p} \cap \pb \right) \cap \surj^{C_p} = \{\id_\emptyset\}^{C_p} \cap \pb = \{\id_\emptyset\}^{C_p}.
	\)
\end{proof}
\begin{lemma}
	\label{(inite cup idemptyCp) cap surje}
	$\left(\init^e \cup \{\id_\emptyset\}^{C_p} \right) \cap \surj^e = \{\id_\emptyset\}^{C_p} \cap \surj^e$
\end{lemma}

\begin{proof}
	We have
	\begin{align*}
		\left(\init^e \cup \{\id_\emptyset\}^{C_p} \right) \cap \surj^e
		 & =
		\left(\init^e \cap \surj^e \right) \cup \left(\{\id_\emptyset\}^{C_p} \cap \surj^e\right)
		\\
		 & =
		\{\id_\emptyset\} \cup \left(\{\id_\emptyset\}^{C_p} \cap \surj^e\right)
		\\
		 & =
		\{\id_\emptyset\}^{C_p} \cap \surj^e,
	\end{align*}
	where the last line follows because $\id_\emptyset$ is in both classes $\{\id_\emptyset\}^{C_p}$ and $\surj^e$.
\end{proof}

\begin{lemma}
	\label{(initCp cap inje) cap surje}
	$\left(\init^{C_p} \cap \inj^e\right) \cap \surj^e = \{\id_\emptyset\}^{C_p} \cap \bij^e$
\end{lemma}
\begin{proof}
	The left-hand side is equal to $\init^{C_p} \cap \bij^e$. A morphism $f : X \to Y$ lies in $\init^{C_p} \cap \bij^e$ if and only if $X^{C_p} = \emptyset$ and $f \in \bij^e$ (i.e. $f$ is an isomorphism). Given that $f$ is an isomorphism, the condition $X^{C_p} = \emptyset$ is equivalent to $Y^{C_p} = \emptyset$. Hence, $f \in \init^{C_p} \cap \bij^e$ if and only if $f \in \{\id_\emptyset\}^{C_p} \cap \bij^e$.
\end{proof}

\begin{lemma}
	\label{(initCp cap surje) cap inje}
	$\left(\init^{C_p} \cap \surj^e\right) \cap \inj^e = \{\id_\emptyset\}^{C_p} \cap \bij^e$
\end{lemma}
\begin{proof}
	Identical to \Cref{(initCp cap inje) cap surje}.
\end{proof}

\begin{lemma}
	\label{(inite cup idemptyCp) cap surjCp}
	$\left(\init^e \cup \{\id_\emptyset\}^{C_p} \right) \cap \surj^{C_p} = \{\id_\emptyset\}^{C_p}$
\end{lemma}

\begin{proof}
	We have
	\begin{align*}
		\left(\init^e \cup \{\id_\emptyset\}^{C_p} \right) \cap \surj^{C_p}
		                                                                               & =
		\left(\init^e \cap \surj^{C_p} \right) \cup \left(\{\id_\emptyset\}^{C_p} \cap \surj^{C_p}\right)
		\\
		                                                                               & =
		\left(\{\id_\emptyset\}^{C_p} \cap \init^e\right) \cup \{\id_\emptyset\}^{C_p} & \text{\cref{inite cap surjCp}}
		\\
		                                                                               & =
		\{\id_\emptyset\}^{C_p}
	\end{align*}
	where the last line follows because $\{\id_\emptyset\}^{C_p} \cap \init^e \subseteq \{\id_\emptyset\}^{C_p}$.
\end{proof}

\begin{lemma}
	\label{(initCp cap inje) cap surjCp}
	$\left(\init^{C_p} \cap \inj^e\right) \cap \surj^{C_p} = \{\id_\emptyset\}^{C_p} \cap \inj^e$
\end{lemma}
\begin{proof}
	Rearranging and applying \cref{initCp cap surjCp}, we have \[\left(\init^{C_p} \cap \inj^e\right) \cap \surj^{C_p} = \left(\init^{C_p} \cap \surj^{C_p}\right) \cap \inj^e = \{\id_\emptyset\}^{C_p} \cap \inj^e.\qedhere\]
\end{proof}

\begin{lemma}
	\label{(initCp cap surje) cap surjCp}
	$\left(\init^{C_p} \cap \surj^e\right) \cap \surj^{C_p} = \{\id_\emptyset\}^{C_p} \cap \surj^e$
\end{lemma}
\begin{proof}
	By \cref{initCp cap surjCp}, we have $\left(\init^{C_p} \cap \surj^e\right) \cap \surj^{C_p} = \{\id_\emptyset\}^{C_p} \cap \surj^e$.
\end{proof}

\begin{lemma}
	\label{(initCp cap pb) cap surje}
	$\left(\init^{C_p} \cap \pb \right) \cap \surj^e = \{\id_\emptyset\}^{C_p} \cap \surj^e$
\end{lemma}
\begin{proof}
	To begin, the right-hand-side $ \{\id_\emptyset\}^{C_p} \cap \surj^e$ is contained in the left because $\{\id_\emptyset\}^{C_p}$ is contained in both $\init^{C_p}$ and $\pb$ by \cref{easier pb and homog for Cp}.

	Conversely, let $f \colon X \to Y$ such that $f \in \init^{C_p} \cap \pb \cap \surj^e$. Then $f^{C_p}$ has domain $\emptyset$, and $f \in \pb$ implies that $f^{-1}(Y^{C_p}) = \emptyset$. Now $f \in \surj^e$ implies that $Y^{C_p} = \emptyset$.  Therefore, $f^{C_p} = \id_\emptyset$. Hence, $f \in \{\id_\emptyset\}^{C_p}$, and $f \in \surj^e$ by assumption.
\end{proof}

\begin{lemma}
	\label{(initCp cap surje) cap pb}
	$\left(\init^{C_p} \cap \surj^e\right) \cap \pb = \{\id_\emptyset\}^{C_p} \cap \surj^e$
\end{lemma}
\begin{proof}
	Identical to \cref{(initCp cap pb) cap surje}.
\end{proof}

\begin{lemma}
	\label{((inite cup idemptyCp) cap inje) cap surje}
	$\left(\left(\init^e \cup \{\id_\emptyset\}^{C_p}\right) \cap \inj^e \right) \cap \surj^e = \{\id_\emptyset\}^{C_p} \cap \bij^e$
\end{lemma}
\begin{proof}
	We have
	\begin{align*}
		\left(\left(\init^e \cup \{\id_\emptyset\}^{C_p}\right) \cap \inj^e \right) \cap \surj^e
		 & =
		\left(\init^e \cup \{\id_\emptyset\}^{C_p}\right) \cap \bij^e
		\\
		 & =
		\left(\init^e \cap \bij^e \right) \cup \left(\{\id_\emptyset\}^{C_p} \cap \bij^e\right)
		\\
		 & =
		\{\id_\emptyset\} \cup \left(\{\id_\emptyset\}^{C_p} \cap \bij^e\right)
		\\
		 & =
		\{\id_\emptyset\}^{C_p} \cap \bij^e,
	\end{align*}
	where the last line follows because $\id_\emptyset$ is a bijection and is in $\{\id_\emptyset\}^{C_p}$.
\end{proof}

\begin{lemma}
	\label{((inite cup idemptyCp) cap inje) cap surjCp}
	$\left(\left(\init^e \cup \{\id_\emptyset\}^{C_p}\right) \cap \inj^e \right) \cap \surj^{C_p} = \{\id_\emptyset\}^{C_p} \cap \inj^e$
\end{lemma}
\begin{proof}
	We have
	\begin{align*}
		\rlap{$\left(\left(\init^e \cup \{\id_\emptyset\}^{C_p}\right) \cap \inj^e \right) \cap \surj^{C_p}$}\phantom{xyz} \\
		 & =
		\left(\init^e \cap \inj^e \cap \surj^{C_p} \right)\cup \left(  \{\id_\emptyset\}^{C_p} \cap \inj^e \cap \surj^{C_p} \right)
		\\
		 & =
		\left(\init^e \cap \surj^{C_p} \right) \cup \left(\{\id_\emptyset\}^{C_p} \cap \inj^e \right)
		\\
		 & =
		\left(\{\id_\emptyset\}^{C_p} \cap \init^e \right) \cup \left(\{\id_\emptyset\}^{C_p} \cap \inj^e \right)
		 &
		\text{\cref{inite cap surjCp}}
		\\
		 & =
		\{\id_\emptyset\}^{C_p} \cap \left(\init^e \cup \inj^e\right)
		\\
		 & =
		\{\id_\emptyset\}^{C_p} \cap \inj^e.
	\end{align*}
	Above, we have freely used the facts that $\init^e \subseteq \inj^e$ and $\{\id_\emptyset\}^{C_p} \subseteq \surj^{C_p}$.
\end{proof}

\vfill

\section{Reference for Classes of Morphisms Used in \Cref{section: Cp-sets}}
\label{fig:cheat sheet}

	\centering
	\small
	\begin{tblr}{width=\linewidth, colspec={Q[l,wd=0.30\linewidth]X[l]}, colsep=4pt, rowsep=3pt}
		\toprule
		\textbf{Name} & \textbf{Definition}                                                                                                                                                                                                                        \\
		\midrule
		Primordial vacuous class
		              & One of $\{\id_\emptyset\}$, $\{\id_\emptyset\}^{C_p}$, $\init^e$, $\init^{C_p}$, or $\init^{C_p}\cap\pb$.                                                                                                                                  \\
		\addlinespace
		Collateral class
		              & A nonempty union of primordial vacuous classes. There are six collateral classes: the five primordial vacuous classes, and $\init^e\cup\{\id_\emptyset\}^{C_p}$.                                                                           \\
		\addlinespace
		Primordial class
		              & One of $\inj^e$, $\inj^{C_p}$, $\surj^e$, $\surj^{C_p}$, $\pb$, or $\homog$.                                                                                                                                                               \\
		\addlinespace
		Primordial intersection
		              & A (possibly empty) intersection of primordial classes.                                                                                                                                                                                     \\
		\addlinespace
		Mixed class
		              & A union $A\cup A'$, where $A$ is a primordial intersection and $A'$ is a collateral class. \textbf{Every right-lifting class $\rlp{\{f\}}$ of a single morphism is a mixed class, but not every mixed class is of the form $\rlp{\{f\}}$.} \\
		\addlinespace
		Collateral intersection
		              & An intersection of one or more collateral classes with a primordial intersection.                                                                                                                                                          \\
		\addlinespace
		Arbitrary right-lifting class
		              & Any right-lifting class of a set of morphisms, $\mathcal{M}^\boxslash = \bigcap_{f \in \mathcal{M}} \rlp{\{f\}}$.
		Each $\rlp{\{f\}}$ has the form $A \cup A'$, and distributing the intersection over these decompositions shows that any right-lifting class has the form
		$A\cup\bigcup_k Q_k$, with $A$ a primordial intersection and each $Q_k$ a collateral intersection.                                                                                                                                                         \\
		\bottomrule
	\end{tblr}

\newpage
\nocite{sagemath}
\bibliographystyle{amsalpha}
\bibliography{sources.bib}

\end{document}